\documentclass[twoside,a4paper,12pt]{article}
\usepackage{amsmath}
\usepackage[usenames,dvipsnames]{color}

\newcommand{\rz}{{\rm I\!R}}
\newcommand{\uad}{{\mathcal{U}_{\rm ad}}}
\newcommand{\nz}{{\rm I\!N}}
\newcommand{\ve}{{\varepsilon}}
\newcommand{\heins}{{H^1(\Omega)}}
\newcommand{\hzwei}{{H^2(\Omega)}}

\newcommand{\heg}{{H^1(0,T;L^2(\Gamma))}}
\newcommand{\lio}{{L^\infty(\Omega)}}
\newcommand{\liq}{{L^\infty(Q)}}
\newcommand{\lig}{{L^\infty(\Gamma)}}
\newcommand{\lis}{{L^\infty(\Sigma)}}
\newcommand{\lzo}{{L^2(\Omega)}}
\newcommand{\lzq}{{L^2(Q)}}
\newcommand{\qed}{\hfill
\colorbox{black}{\hspace{-0.01cm}}}
\renewcommand{\min}{\mathop{\rm Min}}

\newcommand\testopari{\sc Colli \ --- \ Gilardi \ --- \ Sprekels}
\newcommand\testodispari{\sc Optimal Boundary
Control of a Nonstandard Phase Field System}
\def\hth{H^{3/2}(\Omega)}

\begin{document}
\thispagestyle{empty}
\begin{center}
{\bf {\huge Analysis and Optimal Boundary\\[1mm]
Control of a Nonstandard System\\[2mm] 
of Phase Field Equations}\\[6mm]

\vspace{9mm}
{\large\bf Pierluigi Colli$\!^1\!$,
Gianni Gilardi\footnote{Dipartimento di Matematica  ``F. Casorati'',
Universit\`a di Pavia, Via Ferrata, 1, 27100 Pavia,  Italy,
e-mail: pierluigi.colli@unipv.it, gianni.gilardi@unipv.it },\\[2mm]
and 
J\"urgen Sprekels\footnote{Weierstrass Institute for 
Applied Analysis and Stochastics,
Mohrenstr. 39, 10117 Berlin, Germany,
e-mail: juergen.sprekels@wias-berlin.de \\[2mm]
Key words: Nonlinear phase field systems, 
Cahn--Hilliard systems, parabolic systems, optimal boundary control, 
first-order necessary optimality conditions.\\ 
AMS (MOS) Subject Classification: 74A15, 35K55, 49K20.}}
}
\end{center}

\vspace{7mm}
{\bf Abstract. }
We investigate a nonstandard phase field 
model of Cahn-Hilliard type. The model, which was introduced in
\cite{PG}, describes two-species phase segregation and consists of a
system of two highly nonlinearly coupled PDEs. It has been studied
recently in  
\cite{CGPS3}, \cite{CGPS4} for the case of homogeneous Neumann
boundary conditions. In this paper, we investigate the case that the
boundary condition for one of the unknowns of the system is of third
kind and nonhomogeneous. For the resulting system, we show
well-posedness, and we study optimal boundary control
problems. Existence of optimal controls is shown, and the first-order
necessary optimality conditions are derived. Owing to the strong
nonlinear couplings in the PDE system, standard arguments of optimal
control theory do not apply directly, although the control constraints
and the cost functional will be of standard type.   


\section{Introduction}
Let $\Omega\subset\rz^3$ denote an open and bounded domain whose smooth
boundary $\Gamma$ has outward unit normal ${\bf n}$, let $T>0$ be a given final time, and 
let $Q:=\Omega\times (0,T)$,
$\Sigma:=\Gamma\times (0,T)$. In this paper, we study the following initial-boundary value problem:
\begin{eqnarray}
\label{eq:1.1}
(\varepsilon+2\,\rho)\mu_t+\mu\rho_t-\Delta\mu=0\,\quad\mbox{a.\,e. in }\,Q,\\
\label{eq:1.2}
\delta\rho_t-\Delta\rho+f'(\rho)=\mu\,\quad\mbox{a.\,e. in }\,Q,\\
\label{eq:1.3}
\frac{\partial \rho}{\partial \bf n}=0\,,\quad\,\frac{\partial \mu}{\partial \bf n}=
\alpha(u-\mu)\,\quad\mbox{a.\,e. on }\Sigma,\\
\label{eq:1.4}
\rho(x,0)=\rho_0(x)\,,\quad \mu(x,0)=\mu_0(x)\,,\quad \mbox{for a.\,e. } \, x\in  \Omega.
\end{eqnarray}

The PDE system (\ref{eq:1.1})--(\ref{eq:1.2}) constitutes a phase field model of Cahn-Hilliard type 
that  describes  phase segregation of two species (atoms and vacancies, say) on a lattice
in the presence of diffusion. It has been introduced recently in \cite{PG} and  \cite{CGPS3};
for the general physical background, we refer the reader to \cite{PG}. 
The unknown variables  are the  {\em order parameter} $\rho$,
interpreted as a volumetric density,  and the {\em chemical potential} \,$\mu$. For physical reasons, we must 
have $0\le \rho\le 1$ and $\mu > 0$ almost everywhere in~$Q$. The boundary (control)
function $u$ on the right-hand side of (\ref{eq:1.3})$_2$ plays the role of a {\em microenergy source}.  
Moreover, $\varepsilon$ and $\delta$ are positive constants, and the nonlinearity $f$ is a double-well potential defined in~$(0,1)$, whose derivative 
$f'$ is singular at the endpoints $\rho=0$ and $\rho=1$; a typical example is $f=f_1+f_2$, with
$f_2$ smooth and $f_1(\rho)=c\,(\rho\,\log(\rho)+(1-\rho)\,\log(1-\rho))$,  where $\,c\,$ is a positive
constant.

The PDE system (\ref{eq:1.1})--(\ref{eq:1.4}) is singular, with highly nonlinear
and nonstandard coupling. 
In particular, unpleasant nonlinear couplings involving time derivatives occur in (\ref{eq:1.1}),
and the expression $f'(\rho)$ in (\ref{eq:1.2}) may become singular.  
In the recent papers \cite{CGPS3}, \cite{CGPS4}, well-posedness and asymptotic behavior for $t\to
\infty$ and $\varepsilon\searrow 0$ of the system (\ref{eq:1.1})--(\ref{eq:1.4}) were established for the 
case when the second boundary condition in (\ref{eq:1.3}) 
is replaced by the homogeneous Neumann boundary condition $\,\partial\mu/\partial{\bf n}=0$; 
a distributed optimal control problem for this situation was  analyzed in \cite{CGPS5}. We also refer
to the papers \cite{CGPS1} and \cite{CGPS2}, where the corresponding Allen-Cahn model
was discussed.   

The paper is organized as follows: 
in Section~2, we state the general assumptions and prove the existence of a strong solution to the problem. 
Section~3 is concerned with the issues of uniqueness and stability. 
Section~4
then brings the study of a boundary control problem for the system 
(\ref{eq:1.1})--(\ref{eq:1.4}). We show existence of a solution to the
optimal control problem and derive the first-order necessary optimality conditions, as usual given in
terms of the adjoint system and a variational inequality. 

Throughout the paper, we make repeated use of H\"older's inequality, of the elementary Young inequality
\begin{equation} 
\label{eq:1.5}
a\,b\le \gamma a^2 + \frac 1 {4\,\gamma}\,b^2, \quad\mbox{for every }\,a,b\ge 0\,\mbox{ and }\,\gamma>0,
\end{equation}
of the interpolation inequality
\begin{eqnarray} 
\label{eq:1.6}
&&\|v\|_{L^r(\Omega)}\,\le\,\|v\|_{L^p(\Omega)}^\theta\,\|v\|_{L^q(\Omega)}^{1-\theta}
\,\quad\forall \,v\in L^p(\Omega)
\cap L^q(\Omega),\nonumber\\[1mm]
&&\mbox{where }\,p,q,r\in [1,+\infty], \quad \theta\in [0,1], \quad\mbox{and }\,\,\,\frac 1 r =\frac \theta p
\,+\,\frac{1-\theta} q\,\,,\qquad 
\end{eqnarray}
and, since $\mbox{dim}\,\Omega\le 3$, of the continuity of the embeddings $H^1(\Omega)\subset L^q(\Omega)$ for
$1\le q\le 6$, where, with
constants $\hat C_q>0$ depending only on $\Omega$,
\begin{equation}
\label{eq:1.7}
\|v\|_{L^q(\Omega)}\,\le\,\hat C_q\,\|v\|_{\heins}\,\quad\forall\,v\in\heins\,,\quad\,1\le q\le 6,
\end{equation}
and where the embeddings are compact for $1\le q<6$.   
We also use the Sobolev spaces $H^s(\Omega)$ of real order~$s>0$
and recall the compact embeddings 
$H^s(\Omega)\subset H^1(\Omega)$ and $H^s(\Omega)\subset C(\overline\Omega)$
for $s>1$ and $s>3/2$, respectively, and, e.\,g., the estimate,
with a constant $\hat C_\infty>0$ depending only on~$\Omega$,
\begin{equation}
\label{eq:1.8}
\|v\|_{C(\overline \Omega)}\,\le\,\hat C_\infty\,\|v\|_{\hzwei}\,\quad\forall\,v\in\hzwei\,.
\end{equation}

\section{Problem statement and existence}
\setcounter{equation}{0}
Consider the initial-boundary value  problem
(\ref{eq:1.1})--(\ref{eq:1.4}). For convenience, 
we introduce the abbreviated
notation 
$$H=L^2(\Omega), \quad \ V=H^1(\Omega), \quad\ W=\left\{w\in H^2(\Omega)\, :
\ \ \partial w/\partial {\bf n}=0 \ \mbox{on }\,\Gamma\right\}. 
$$
We endow these spaces with their 
standard norms, for which we use self-explaining notation like $\|\cdot\|_V$; for simplicity, we
also write $\|\cdot\|_H$ for the norm in the space $H\times H\times H$. Recall that the embeddings
$W\subset V\subset H$ are compact. Moreover, since $V$ is dense in $H$, we can identify $H$
with a subspace of $V^*$ in the usual way, i.\,e., by setting $\langle u,v\rangle_{V^*,V}
=(u,v)_H$ for all $u\in H$ and $v\in V$, 
where $\langle \cdot\,,\,\cdot\rangle_{V^*,V}$ denotes the duality pairing
between $V^*$ and $V$.
Then also the embedding $H\subset V^*$ is compact. 

\vspace{2mm}
We make the following assumptions on the data:

\vspace{3mm}
(A1) $f=f_1+f_2$, where $f_1\in C^2(0,1)$ is convex, $f_2\in C^2[0,1]$,  and
\begin{equation}
\label{eq:2.1}
\lim_{r\searrow 0} f_1'(r)=-\infty, \quad \lim_{r\nearrow 1} f_1'(r)= +\infty.
\end{equation}
(A2) $\rho_0\in W$, $f'(\rho_0) \in H$, $\mu_0\in V$, and
\begin{equation}
\label{eq:2.2}
0<\rho_0(x)<1 \quad \forall\,x\in\Omega,\,  \quad \mu_0\ge 0\,\,
\mbox{ a.\,e. in }\,\Omega. 
\end{equation}
(A3) $u\in H^1(0,T;L^2(\Gamma))$, and $u\ge 0$ \,a.\,e. on $\,\Sigma$.

\vspace{2mm}
(A4) $\alpha\in L^\infty(\Gamma)$, and $\,\alpha(x)\ge\alpha_0>0$ \,for 
almost every $\,x\in \Gamma$. 

\vspace{2mm}
Notice that (A2) implies that $\rho_0\in C(\overline\Omega)$ and, thanks to 
the convexity of $\,f_1\,$, also that $f(\rho_0)\in H$. 

The following existence result resembles that of Theorem~2.1 in \cite{CGPS3}.

\vspace{5mm}
{\bf Theorem 2.1} \quad {\em Suppose that the hypotheses} (A1)--(A4) {\em are 
satisfied. Then the system} (\ref{eq:1.1})--(\ref{eq:1.4}) {\em has a solution}
$(\rho,\mu)$ {\em such that}
\begin{eqnarray}
\label{eq:2.3}
&&\rho\in W^{1,\infty}(0,T;H)\cap H^1(0,T;V)\cap L^\infty(0,T;W),\\[1mm]
\label{eq:2.4}
&&\mu\in H^1(0,T;H)\cap C^0([0,T];V)\cap L^2(0,T;\hth),\\[1mm]
\label{eq:2.5}
&& f'(\rho)\in L^\infty(0,T;H),\\[1mm]
\label{eq:2.6}
&&0<\rho<1 \quad a.\,e.\,\, in \,\,Q, \quad\,\,\mu\ge 0\quad a.\,e.\,\, in \,\,Q.
\end{eqnarray} 

\vspace{3mm}
{\bf Remark 2.2} \quad 
The $H^{3/2}$ space regularity for $\mu$ 
is optimal due to the $L^2$ space regularity of~$u$ given by~(A3).
Nevertheless, both equation (\ref{eq:1.1}) and the boundary condition 
for $\mu$ contained in (\ref{eq:1.3}) can be understood a.e.\ in~$Q$
and a.e.\ on~$\Sigma$, respectively,
and the standard integration by parts is correct, as we briefly explain
(so~that we can both refer to that formulation and use integration by parts).
In principle, one can replace the equation and the boundary condition
by the usual variational formulation,
namely
\begin{equation*}
\int_\Omega \left[ (\varepsilon+2\,\rho)\mu_t+\mu\rho_t\right] \, v \,dx
+\int_\Omega \nabla\mu\cdot\nabla v \,dx
+\int_\Gamma \alpha(\mu-u) v \,d\sigma
=0
\end{equation*}
(where $d\sigma$ stands for the surface measure)
for every $v\in V$, a.e.\ in $(0,T)$,
or an integrated-in-time version of it.
This implies that (\ref{eq:1.1}) is satisfied in the sense of distributions,
whence $\Delta\mu$ belongs to $L^2(Q)$ by comparison,
and the equation can be understood~a.e.\ in~$Q$, a~posteriori.
The last regularity (\ref{eq:2.4}) of $\mu$ 
and the condition $\Delta\mu\in L^2(0,T;L^2(\Omega))$ just observed 
also ensure that the trace
$\frac{\partial\mu}{\partial\bf n}|_\Sigma$ has a meaning in the space 
$L^2(0,T;L^2(\Gamma))$
due to the trace theorem \cite[Thm.~7.3]{LioMag}
(we~just observe that the space $\Xi^{-1/2}(\Omega)$ 
that enters such a result is larger than~$L^2(\Omega)$),
so that the boundary condition can be read a.e.\ on~$\Sigma$.

\vspace{3mm}
{\em Proof of Theorem 2.1.} 
The proof follows closely the lines of the proof of Theorem~2.1 in \cite{CGPS3}, where a homogeneous Neumann boundary condition for $\mu$ was investigated.

\vspace{2mm}
\underline{Step 1: Approximation.} \,\,We employ an approximation
scheme based on a time delay in the right-hand side of
(\ref{eq:1.2}). To this end, we introduce for $\tau>0$ the translation
operator $\,\mathcal{T}_\tau:L^1(0,T;H)\to L^1(0,T;H)$, which for $v\in L^1(0,T;H)$ and almost every $t\in (0,T)$ is defined by
\begin{equation}
\label{eq:2.7}
( \mathcal{T}_\tau)(t):=v(t-\tau) \quad\mbox{if }\, t>\tau, \quad\mbox{and } \,\,(\mathcal{T}_\tau)(t):=\mu_0
\quad\mbox{if }\, t\le\tau\,.  
\end{equation}
Now, let $N\in\nz$ be arbitrary, and $\tau:=T/N$. We seek functions $(\rho^\tau,\mu^\tau)$ satisfying (\ref{eq:2.3})--(\ref{eq:2.6}) (with $(\rho,\mu)$ replaced by $(\rho^\tau,\mu^\tau)$),
which solve the system 
\begin{eqnarray}
\label{eq:2.8}
(\varepsilon+2\,\rho^\tau)\mu_t^\tau+\mu^\tau\rho^\tau_t-\Delta\mu^\tau=0\,\quad\mbox{a.\,e. in }\,Q,\\[1mm]
\label{eq:2.9}
\delta\rho^\tau_t-\Delta\rho^\tau+f'(\rho^\tau)=\mathcal{T}_\tau\mu^\tau\,\quad\mbox{a.\,e. in }\,Q,\\[1mm]
\label{eq:2.10}
\frac{\partial \rho^\tau}{\partial \bf n}=0\,,\quad\,\frac{\partial \mu^\tau}{\partial \bf n}=
\alpha(u-\mu^\tau)\,\quad\mbox{a.\,e. on }\Sigma,\\[1mm]
\label{eq:2.11}
\rho^\tau(x,0)=\rho_0(x)\,,\quad \mu^\tau(x,0)=\mu_0(x)\,,\quad\mbox{for a.\,e. } \, x\in  \Omega.
\end{eqnarray}

We note that Remark~2.2 also applies to the approximating problem.
To prove the existence of a solution, we put $t_n:=n\,\tau$, $I_n:=[0,t_n]$, $1\le n\le N$, and consider for $1\le n\le N$ the problem
 \begin{eqnarray}
\label{eq:2.12}
(\varepsilon+2\,\rho^n)\mu_t^n+\mu^n\rho^n_t-\Delta\mu^n=0\,\quad\mbox{a.\,e. in }\,\Omega\times I_n,\\[1mm]
\label{eq:2.13}
\mu^n(0)=\mu_0\quad\mbox{a.\,e. in } \,\Omega\,, \quad \frac{\partial \mu^n}{\partial \bf n}=
\alpha(u-\mu^n)\,\quad\mbox{a.\,e. on }\Gamma\times I_n,\\[1mm]
\label{eq:2.14}
\delta\rho^n_t-\Delta\rho^n+f'(\rho^n)=\mathcal{T}_\tau\mu^{n-1}\,\quad\mbox{a.\,e. in }\,\Omega\times I_n,\\[1mm]
\label{eq:2.15}
\rho^n(0)=\rho_0\quad\mbox{a.\,e. in } \,\Omega\,,\quad \frac{\partial \rho^n}{\partial \bf n}=0\,,\quad\mbox{a.\,e. on }\Gamma\times I_n\,.
\end{eqnarray}

Notice that the operator $\mathcal{T}_\tau$ acts on functions that are not defined on the entire interval
$(0,T)$. However, its meaning is still given by (\ref{eq:2.7}) if $n>1$, and for $n=1$ we simply put $\mathcal{T}_\tau\mu^{n-1}=\mu_0$. 

Clearly, we have $(\rho^\tau,\mu^\tau)=(\rho^N,\mu^N)$ if $(\rho^N,\mu^N)$ exists. We claim that the systems (\ref{eq:2.12})--(\ref{eq:2.15}) can be uniquely solved by induction for $n=1,...,N$, where, for $1\le n\le N$,
\begin{eqnarray}
\label{eq:2.16}
&&\rho^n\in W^{1,\infty}(I_n;H)\cap H^1(I_n;V)\cap L^\infty(I_n;W),\\[1mm]
\label{eq:2.17}
&&\mu^n\in H^1(I_n;H)\cap C^0(I_n;V)\cap L^2(I_n;\hth),\\[1mm]
\label{eq:2.18} 
&&0<\rho^n<1 \quad\mbox{a.\,e. in }\,\Omega\times I_n, \quad\,\,\mu^n \ge 0 \quad \mbox{a.\,e. in }\,\Omega\times I_n.
\end{eqnarray} 

To prove the claim, suppose that for some $\,n\in\{1,\ldots,N\}\,$ the problem (\ref{eq:2.12})--(\ref{eq:2.15}) 
has a unique solution satisfying (\ref{eq:2.16})--(\ref{eq:2.18}), 
where the index $n$ is replaced by $n-1$. 
Then it follows with exactly the same argument as in the proof of Theorem~2.1 in \cite{CGPS3} that 
the initial-boundary value problem (\ref{eq:2.14}), (\ref{eq:2.15}) has  a unique solution $\rho^n$ 
that satisfies (\ref{eq:2.16}) and the first inequality in (\ref{eq:2.18}). 
Substituting $\rho^n$ in  (\ref{eq:2.12}), we infer that the linear initial-boundary value problem (\ref{eq:2.12}), (\ref{eq:2.13}) 
has a unique solution $\mu^n$ satisfying (\ref{eq:2.17}). 
Notice here that the regularity of $\mu^n$ 
follows from the fact that $u\in \heg$.

It remains to show that $\mu^n$ is nonnegative almost everywhere. To this end, we test (\ref{eq:2.12}) by 
$-(\mu^n)^-$, where $(\mu^n)^-$ denotes the negative part of $\mu^n$. 
Using integration by parts and the boundary condition in (\ref{eq:2.13}), we obtain the identity
\begin{eqnarray*}
&&\frac 1 2 \int_0^t\!\!\int_\Omega \frac d {dt} \Bigl((\ve+2\rho^n)\,\left|(\mu^n)^-\right|^2 \Bigr)\,dx\,ds
+ \int_0^t\!\!\int_\Omega\left|\nabla(\mu^n)^-\right|^2\,dx\,ds\\[1mm]
&&+\,\int_0^t\!\!\int_\Gamma \alpha\,\left|(\mu^n)^-\right|^2 \,d\sigma\,ds\,+\,
\int_0^t\!\!\int_\Gamma \alpha\,u\,(\mu^n)^-\,d\sigma\,ds\,=\,0\,.
\end{eqnarray*}

From the fact that $\rho^n$, $\rho_0$, $\mu_0$, $\alpha$, $u$\, are all nonnegative, we infer that 
\begin{eqnarray*}
&&\ve \,\int_\Omega \left|(\mu^n)^-(t)\right|^2 \,dx \,\le\, \int_\Omega  
(\ve+2\rho^n(t))\,\left|(\mu^n)^-(t)\right|^2 \,dx\\[1mm]
&&\le\, \int_\Omega (\ve+2\rho_0)\left|\mu_0^-\right|^2\,dx =\,0\,.
\end{eqnarray*}

Hence, $(\mu^n)^-=0$, i.\,e., $\mu^n\ge 0$ a.\,e. in $\Omega\times I_n$, and the claim is proved.

\vspace{5mm}
\underline{Step 2: A priori estimates.} \,\,Now that the well-posedness of the problem \linebreak (\ref{eq:2.8})--(\ref{eq:2.11})
is established, we perform a number of a priori estimates for its solution. For the sake of a better readability,
we will omit the index $\tau$ in the calculations. In what follows, we denote by $C>0$ positive constants that may depend 
on the data of the system but not on $\tau$. The meaning of $C$ may change from line to line and even in the same chain
of inequalities.

\vspace{3mm}
\underline{First estimate.}  \,\,Since $\,\partial_t\bigl((\ve/2)\mu^2+\rho\mu^2\bigr)=\bigl((\ve+2\rho)\mu_t
+\mu\rho_t\bigr)\,\mu$, testing of (\ref{eq:2.8}) by $\mu$ yields, for every $t\in [0,T]$,
\begin{eqnarray*}
&&\int_\Omega\Bigl(\frac \ve 2 \,\mu^2+\rho\mu^2\Bigr)(t)\,dx \,+\,\int_0^t\!\!\int_\Omega |\nabla\mu|^2\,dx\,ds
\,+\,\int_0^t\!\!\int_\Gamma\alpha\,\mu^2\,d\sigma\,ds\nonumber\\[1mm]
&&=\,\int_\Omega\Bigl(\frac \ve 2 \,\mu_0^2+\rho_0\mu_0^2\Bigr)(t)\,dx  \,+\,\int_0^t\!\!\int_\Gamma\alpha\,u\,\mu\,d\sigma\,ds,
\end{eqnarray*} 

whence, using Young's inequality and (A2)--(A4), we can conclude that
\begin{equation}
\label{eq:2.19}
\|\mu\|_{L^\infty(0,T;H)\cap L^2(0,T;V)}\,\le\,C\,.
\end{equation}

\vspace{3mm}
\underline{Second estimate.}  \,\,Next, we test (\ref{eq:2.9}) by $\rho_t$. Applying (\ref{eq:2.19}), recalling the fact that $f(\rho_0)\in H$, and invoking Young's inequality, we easily see that

\begin{equation}
\label{eq:2.20}
\|\rho\|_{H^1(0,T;H)\cap L^\infty(0,T;V)}\,+\,\|f(\rho)\|_{L^\infty(0,T;L^1(\Omega))}\,\le\,C\,.
\end{equation}

\vspace{5mm}
\underline{Third estimate.}  \,\,We rewrite Eq. (\ref{eq:2.9}) in the form
\begin{equation*}
-\Delta\rho\,+\,f_1'(\rho)\,=\,-\,\delta\,\rho_t-f_2'(\rho)\,+\,\mathcal{T}_\tau\mu\,
\end{equation*}

and observe that the right-hand side is bounded in $L^2(Q)$. Hence, applying a standard procedure
(e.\,g., testing by $f_1'(\rho)$), and invoking elliptic regularity, we find that
\begin{equation}
\label{eq:2.21}
\|\rho\|_{L^2(0,T;W)}\,+\,\|f_1'(\rho)\|_{L^2(Q)}\,\le\,C\,.
\end{equation}

\vspace{5mm}
\underline{Fourth estimate.} \,\,We differentiate Eq. (\ref{eq:2.9}) formally with respect to $t$ and test the
resulting equation with $\rho_t$ (this argument can be made rigorous, see \cite{CGPS3}). Since, owing to the convexity of
$f_1$,  $f_1''(\rho)$ is nonnegative almost everywhere, we find the estimate  

\begin{eqnarray}
\label{eq:2.22}
&&\frac \delta 2 \,\|\rho_t(t)\|^2\,+\,\int_0^t\!\!\int_\Omega \left|\nabla\rho_t\right|^2\,dx\,ds\,\le\,
\frac \delta 2\,\|\Delta\rho_0-f_1'(\rho_0)+\mu_0\|^2_H\nonumber\\[1mm]
&&\quad +\max_{0\le\rho\le 1}\,\left|f_2''(\rho)\right|\int_0^t\!\!\int_\Omega\left|\rho_t\right|^2\,dx\,ds\,+\,
\int_0^t\!\!\int_\Omega\left(\partial_t\mathcal{T}_\tau\mu\right)\,\rho_t\,dx\,ds\nonumber\\[1mm]
&&\le \,C\,+\,\int_0^{t-\tau}\!\!\int_\Omega\mu_t(s)\,\rho_t(s+\tau)\,dx\,ds\,.
\end{eqnarray}

In order to estimate the last integral, we substitute for $\mu_t$, using Eq. (\ref{eq:2.8}). It follows, using integration by parts:

\begin{eqnarray}
\label{eq:2.23}
&&\int_0^{t-\tau}\!\!\int_\Omega\mu_t\,\rho_t(\cdot+\tau)\,dx\,ds \,=\,\int_0^{t-\tau}\!\!\int_\Omega
\frac 1 {\ve+2\rho}\,(\Delta\mu -\mu\,\rho_t)\,\rho_t(\cdot+\tau)\,dx\,ds\nonumber\\[1mm]
&&=\int_0^{t-\tau}\!\!\int_\Omega \left[-\,\frac {\nabla\mu} {\ve+2\rho}\,\cdot \nabla \rho_t(\cdot+\tau)
\,+\,\frac {2\rho_t(\cdot+\tau)}{(\ve+2\rho)^2}\,\nabla\mu\,\cdot\,\nabla\rho\right.\nonumber\\[1mm]
&&\qquad\qquad \quad\left. -\,\frac 1 {\ve+2\rho}\,\rho_t\,\mu\,\rho_t(\cdot+\tau)\right]\,dx\,ds\nonumber\\[1mm]
&&\quad -\,\int_0^{t-\tau}\!\!\int_\Gamma
\frac \alpha {\ve+2\rho}\,(u -\mu)\,\rho_t(\cdot+\tau)\,d\sigma\,ds\,.
\end{eqnarray}

Exactly as in the proof of Theorem~2.1 in \cite{CGPS3}, the domain integral in the second and third lines of (\ref{eq:2.23})
can be estimated from above by an expression of the form
\begin{equation}
\label{eq:2.24}
\frac 1 2 \int_0^t\!\!\int_\Omega \left|\nabla\rho_t\right|^2\,dx\,ds \,+\,C\,\Bigl(1\,+\,\int_0^t \|\mu(s)\|_V^2\,
\|\rho_t(s)\|_H^2\,dx\,ds\Bigr)\,.
\end{equation}
Observe that, owing to the inequality (\ref{eq:2.19}), the mapping $\,s\mapsto \|\mu(s)\|_V^2\,$ belongs to $L^1(0,T)$.

Finally, we estimate the boundary term in the last line of Eq. (\ref{eq:2.23}). To this end, recall that by the trace theorem there is a constant $\,c_\Omega>0\,$, independent of $\tau$, such that 
$\,\|v\|_{L^2(\Gamma)}\,\le\,
c_\Omega\,\|v\|_V\,$ for all $v\in V$. Moreover, we have $\rho\ge 0$ and $\,\alpha\in \lig$. Therefore, we obtain that

\begin{eqnarray}
\label{eq:2.25}
&&\Bigl|\int_0^{t-\tau}\!\!\int_\Gamma
\frac \alpha {\ve+2\rho}\,(u -\mu)\,\rho_t(\cdot+\tau)\,d\sigma\,ds\Bigr|\nonumber\\[1mm]
&&\le\,C\int_0^{t-\tau} \|\rho_t(s+\tau)\|_{L^2(\Gamma)}\,\left(\|u(s)\|_{L^2(\Gamma)}\,+\,
\|\mu(s)\|_{L^2(\Gamma)}\right)\,ds\nonumber\\[1mm]
&&\le\,C\int_0^{t-\tau} \|\rho_t(s+\tau)\|_V\,\left(\|u(s)\|_{L^2(\Gamma)}\,+\,
\|\mu(s)\|_V\right)\,ds\nonumber\\[1mm]
&&\le\,\frac 1 4 \int_0^t \|\rho_t(s)\|_V^2\,ds\,+\,C\,.
\end{eqnarray}

Now we may combine the estimates (\ref{eq:2.22})--(\ref{eq:2.25}) and employ Gronwall's inequality to
conclude that

\begin{equation}
\label{eq:2.26}
\|\rho_t\|_{L^\infty(0,T;H)\cap L^2(0,T;V)}\,\le\,C\,.
\end{equation}

The same argument as in the derivation of (\ref{eq:2.22}) then shows that also

 \begin{equation}
\label{eq:2.27}
\|\rho_t\|_{L^\infty(0,T;W)} \,+\, \|f_1'(\rho)\|_{L^\infty(0,T;H)}\,\le\,C\,.
\end{equation}

\vspace{5mm}
\underline{Fifth estimate.}  \,\,We test equation (\ref{eq:2.8}) by $\mu_t$. Formal integration by parts
(this can be made rigorous), using (A3), the trace theorem and Young's inequality, yields:

\begin{eqnarray}
\label{eq:2.28}
&&\ve\int_0^t\!\!\int_\Omega\left|\mu_t\right|^2\,dx\,ds \,+\,\frac 1 2\,\|\nabla\mu(t)\|_H^2\,
\,+\,\int_\Gamma \frac \alpha 2\,|\mu(t)|^2\,d\sigma\nonumber\\[1mm]
&&\le\,C\,+\,\int_0^t\!\!\int_\Gamma \alpha \,u\,\mu_t\,d\sigma\,+\,\int_0^t\!\!\int_\Omega\left|
\mu\,\rho_t\,\mu_t\right|\,dx\,ds\nonumber\\[1mm]
&&\le\,C\,+\,\int_\Gamma \alpha\, u(t)\,\mu(t)\,d\sigma\,-\,\int_0^t\!\!\int_\Gamma \alpha \,u_t\,\mu\,d\sigma\,ds\,+\,\int_0^t\!\!\int_\Omega\left|
\mu\,\rho_t\,\mu_t\right|\,dx\,ds\nonumber\\[1mm]
&&\le\,\frac C \gamma \,+\,\gamma \,\|\mu(t)\|_V^2\,
+\,\int_0^t \|\mu(s)\|_V^2\,ds\,
+\,\int_0^t\!\!\int_\Omega|\mu|\,|\rho_t|\,|\mu_t|\,dx\,ds\nonumber\\[1mm]
&&\le\,\frac C \gamma +\gamma\|\mu(t)\|_V^2\,+\,\int_0^t \|\mu(s)\|_V^2\,ds\,
+\frac \ve 2\!\int_0^t\!\!\|\mu_t(s)\|_H^2\,ds \nonumber \\[1mm]
&&\hskip1cm {}+\,C\!\!\int_0^t\|\rho_t(s)\|^2_{L^4(\Omega)}\,\|\mu(s)\|^2_{L^4(\Omega)}\,ds\nonumber\\[1mm]
&&\le\,\frac C \gamma +\gamma\|\mu(t)\|_V^2\,+\frac \ve 2\!\int_0^t\!\!\|\mu_t(s)\|_H^2\,ds 
\nonumber \\[1mm]
&&\hskip1cm {}\,+\,C\!\!\int_0^t\left( 1+ \|\rho_t(s)\|^2_V\right)\, \|\mu(s)\|^2_V\,ds.\quad
\end{eqnarray}

Hence, using (\ref{eq:2.26}), choosing $\gamma>0$ sufficiently small, and invoking Gronwall's lemma,
we can conclude that

\begin{equation}
\label{eq:2.29}
\|\mu\|_{H^1(0,T;H)\cap L^\infty(0,T;V)}\,\le\,C\,.
\end{equation}

\vspace{5mm}
\underline {Sixth estimate.} \,\, Since $\,0<\rho<1\,$ a.\,e. in $Q$, and using (\ref{eq:2.26}),
(\ref{eq:2.29}) and the continuity of the
embedding $V\subset L^4(\Omega)$, we can estimate as follows:
\begin{eqnarray}
\label{eq:2.30}
\hskip-1cm&&\|(\ve+2\rho)\mu_t+\mu\rho_t\|_{L^2(Q)}\,\le\,C\,\|\mu_t\|_{L^2(Q)}\,+\,\|\mu\|_{L^\infty(0,T;L^4
(\Omega))}\,\|\rho_t\|_{L^2(0,T;L^4(\Omega))}\nonumber\\[2mm]
\hskip-1cm&&\quad\le\,C\,\left(\|\mu_t\|_{L^2(Q)}\,+\,\|\mu\|_{L^\infty(0,T;V)}\,\|\rho_t\|_{L^2(0,T;V)}\right)\,\le\,C\,.
\end{eqnarray}

Comparison in (\ref{eq:2.8}) then shows the boundedness of $\,\Delta\mu\,$ in $\,L^2(Q)$, and it follows from (\ref{eq:2.8}), (A3) and standard elliptic estimates that also

\begin{equation}
\label{eq:2.31}
\|\mu\|_{L^2(0,T;\hth)}\,\le\,C\,.
\end{equation}

\vspace{5mm}
\underline{Step 3: Conclusion of the proof.} \,\,Collecting all the above estimates, it turns out that there is some sequence $\,\tau_k\searrow 0$ such that
\begin{eqnarray*}
\mu^{\tau_k}\to\mu&&\mbox{ weakly star in }\\[1mm]
&& H^1(0,T;H)\cap L^\infty(0,T;V)\cap L^2(0,T;\hth)\,,\\[1mm]
\rho^{\tau_k}\to\rho&&\mbox{weakly star in }\,W^{1,\infty}(0,T;H)\cap H^1(0,T;V)\cap
L^\infty(0,T;W)\,,\\[1mm]
f_1'(\rho^{\tau_k})\to\xi&&\mbox{weakly star in }\, L^\infty(0,T;H)\,.
\end{eqnarray*}

Thanks to the Aubin-Lions lemma (cf., \cite[Thm.~5.1, p.~58]{Aubin}) and similar results 
to be found in \cite[Sect.~8, Cor.~4]{Simon}, 
we also deduce (recall that even $\hth$ is compactly embedded into~$V$) the~strong convergences
\begin{eqnarray*}
\mu^{\tau_k}\to\mu&&\mbox{strongly in }\,C^0([0,T];H)\cap L^2(0,T;V)\,,\\[1mm]
\rho^{\tau_k}\to\rho&&\mbox{strongly in }\,C^0([0,T];V)
\end{eqnarray*}
and the Cauchy conditions (\ref{eq:1.4}) as a consequence.
In particular, employing a standard monotonicity argument (cf., e.\,g., \cite[Lemma 1.3, p.~42]{Barbu}), we conclude
that $\,0<\rho<1\,$ and $\,\xi=f_1'(\rho)\,$ a.\,e. in $Q$. The strong convergence shown above also entails
that $\,f_2'(\rho^{\tau_k})\to f_2'(\rho)\,$ strongly in $C^0([0,T];H)$ (because $\,f_2'\,$ is Lipschitz continuous), and that $\,\mathcal{T}_{\tau_k}\mu^{\tau_k}\to\mu\,$ strongly in~$\,L^2(Q)$.

Now notice that the above convergences imply, in particular, that
\begin{eqnarray*}
\rho^{\tau_k}\to \rho &&\mbox{strongly in
}\,C^0([0,T];L^6(\Omega))\,,\\[1mm]
\rho^{\tau_k}_t\to \rho_t&&\mbox{weakly in }\,L^2(0,T;L^4(\Omega)),\\[1mm]
\mu^{\tau_k}\to \mu&&\mbox{strongly in }\,
L^2(0,T;L^4(\Omega))\,,\\[1mm]
\mu^{\tau_k}_t\to \mu_t&&\mbox{weakly in }\,L^2(Q)\,.
\end{eqnarray*}
From this, it is easily verified that
\begin{eqnarray*}
\mu^{\tau_k}\,\rho^{\tau_k}_t\to \mu\,\rho_t &&\mbox{weakly in }\,
L^1(0,T;H), \\[1mm]
\rho^{\tau_k}\,\mu^{\tau_k}_t\to\rho\,\mu_t &&\mbox{weakly in }\,
L^2(0,T;L^{3/2}(\Omega)).
\end{eqnarray*}
Now, we are ready to take the limit as $k\to\infty$ in (\ref{eq:2.8})--(\ref{eq:2.10}) 
(written for $\tau=\tau_k$).
Precisely, we can do that as far as $\rho$ is concerned,
while it is easier to take the limit in the variational formulation of (\ref{eq:2.8})
that accounts for the boundary condition (the same as mentioned in Remark~2.2),
or in the following integrated-in-time version of~it
\begin{eqnarray*}
&& \int_0^T\!\!\int_\Omega \left[ (\varepsilon+2\,\rho^\tau)\mu_t^\tau+\mu^\tau\rho^\tau_t\right] \, v \,dx\,dt
+\int_0^T\!\!\int_\Omega \nabla\mu^\tau\cdot\nabla v \,dx\,dt
\\[1mm]
&& \quad +\int_0^T\!\!\int_\Gamma \alpha(\mu^\tau-u) v \,d\sigma\,dt
=0
\qquad
\hbox{for every $v\in L^\infty(0,T;V)$}.
\end{eqnarray*}
Then, we obtain the analogue for~$\mu$, which implies (\ref{eq:1.1}) and (\ref{eq:1.3})$_2$. \qed

\section{Boundedness, uniqueness, and stability}
\setcounter{equation}{0}
In this section, we derive results concerning boundedness, uniqueness and stability of the solutions to system 
(\ref{eq:1.1})--(\ref{eq:1.4}). With respect to boundedness, we have the following result, which resembles Theorem~2.3 in \cite{CGPS3}.

\vspace{5mm}
{\bf Theorem 3.1} \quad {\em Suppose that} (A1)--(A4) {\em are fulfilled, and suppose that the following conditions are satisfied:}

\vspace{2mm}
(A5) \quad$\mu_0\in \lio,\quad  \displaystyle{\inf_{x\in\Omega}\,\rho_0(x)>0, \quad
\sup_{x\in\Omega}\,\rho_0(x)<1.}$

(A6) \quad $u\in\lis$.

\vspace{2mm}
{\em Then any solution $(\rho,\mu)$ of} (\ref{eq:1.1})--(\ref{eq:1.4}) {\em fulfilling} (\ref{eq:2.3})--(\ref{eq:2.6}) {\em also satisfies} 
\begin{equation}
\label{eq:3.1}
\mu\leq\mu^*,
\quad \rho\geq\rho_* \,, 
\quad \hbox{and} \quad
\rho\leq\rho^*
\quad \hbox{a.e.\ in $Q$}
\end{equation} 
{\em for some constants $\mu^*>0$ and $\rho_*\,,\rho^*\in(0,1)$
that depend on the structure of the system and~$T$, on the initial data,
and on an upper bound for the $L^\infty$ norm of~$u$, only.}

\vspace{3mm}
{\em Proof.} \quad Let us just show the boundedness of $\mu$ 
and the first estimate (\ref{eq:3.1}); 
the results for $\rho$ then follow in exactly the same manner as in the proof of Theorem~2.3 in \cite{CGPS3}. Also the result for $\mu$ follows -- up to
some changes that are necessary due to the different boundary condition for $\mu$ -- 
by the same chain of arguments as  in the proof of Theorem~2.3 in \cite{CGPS3}; 
but since this proof does not seem to be standard, we provide it for the reader's convenience. 
So let $(\rho,\mu)$ be any solution to the system (\ref{eq:1.1})--(\ref{eq:1.4}),  (\ref{eq:2.3})--(\ref{eq:2.6}). We set
$$\Phi_0\,:=\,\max\,\{1,\|\mu_0\|_{\lio}\,,\|u\|_{\lis}\}\,,$$
choose any $\,k\in\rz\,$ such that $\,k\ge\Phi_0$, and introduce the auxiliary function
$\,\chi_k\in\liq\,$ by putting, for almost every $(x,t)\in Q$,
$$\chi_k(x,t)=1\quad\mbox{if }\,\mu(x,t)>k,\quad\,\,\mbox{and }\,\,\chi_k(x,t)=0 \quad\mbox{otherwise}.$$ 
Then, we test (\ref{eq:1.1}) by $\,(\mu-k)^+$. We obtain, for any $t\in [0,T]$,
\begin{eqnarray*}
&&\int_\Omega\left(\frac \ve 2+\rho(t)\right)|(\mu(t)-k)^+|^2\,+\,\int_0^t\!\!\int_\Omega|\nabla(\mu-k)^+|^2
\,dx\,ds\\[1mm]
&&\quad+\,\int_0^t\!\!\int_\Gamma \alpha\,(\mu-u)\,(\mu-k)^+\,d\sigma\,ds\\[1mm]
&&=\int_0^t\!\!\int_\Omega\rho_t\,|(\mu-k)^+|^2 \,dx\,ds\,-\,\int_0^t\!\!\int_\Omega\rho_t\,\mu\,(\mu-k)^+
\,dx\,ds\\[1mm]
&&=\,-\,k\int_0^t\!\!\int_\Omega\rho_t\,(\mu-k)^+ \,dx\,ds\,.
\end{eqnarray*}
Now observe that $\alpha$ and $\rho$ are nonnegative and that, by definition of $k$,
$$\alpha\,(\mu-u)\,(\mu-k)^+\,=\,\alpha\,\bigl(|(\mu-k)^+|^2\,+\,(k-u)\,(\mu-k)^+\bigr)\,\ge\,0\quad\mbox{a.\,e. in }\,Q\,.$$
Hence,  using H\"older's inequality, we obtain from the above equality the estimate
\begin{eqnarray*}
&&\frac \ve 2\,\|(\mu(t)-k)^+\|^2_H\,+\,\int_0^t\!\!\int_\Omega|\nabla(\mu-k)^+|^2
\,dx\,ds\\[1mm]
&&\le\,k\int_0^t\|\chi_k(s)\|_{L^{7/2}(\Omega)}\,\|\rho_t(s)\|_{L^{14/3}(\Omega)}\,
\|(\mu-k)^+(s)\|_{\lzo}\,ds\,,
\end{eqnarray*}
whence, using the Gronwall-Bellman lemma as in \cite[Lemma A.4, p. 156]{Brezis},
\begin{eqnarray}
\label{eq:3.2}
&& \Bigl(\ve \,\|(\mu-k)^+\|_{C^0([0,T];H)}^2\,+\,\int_0^T\!\!\int_\Omega|\nabla(\mu-k)^+|^2
\,dx\,dt\Bigr)^{1/2}\nonumber\\[1mm]
&&\le\,\frac k {\sqrt{\ve}}\,\int_0^T \|\chi_k(t)\|_{L^{7/2}(\Omega)}\,\|\rho_t(t)\|_{L^{14/3}(\Omega)}\,dt\nonumber\\[1mm] 
&&\le\,\frac k {\sqrt{\ve}}\,\|\rho_t\|_{L^{7/3}(0,T;L^{14/3}(\Omega))}\,\|\chi_k\|_{L^{7/4}(0,T;
L^{7/2}(\Omega))}\,.
\end{eqnarray}
Next, we apply the continuity of the embedding $V\subset L^6(\Omega)$ 
and the interpolation inequality (\ref{eq:1.6}) with $p=2$, $q=6$, $r=14/3$, and $\theta=1/7$. 
It follows that
\begin{eqnarray*}
&&\|\rho_t\|_{L^{7/3}(0,T;L^{14/3}(\Omega))}\,\le\,\Bigl(\int_0^T\|\rho_t(t)\|_{L^2(\Omega)}^{1/3}\,\|\rho_t(t)\|_{L^6(\Omega)}^2\,dt\Bigr)^{3/7}\nonumber\\[1mm]
&&\le \|\rho_t\|_{L^\infty(0,T;H)}^{1/7}\,\Bigl(\int_0^T\|\rho_t(t)\|_{L^6(\Omega)}^2\,dt\Bigr)^{3/7}
\,\le\,C\,\|\rho_t\|_{L^2(0,T;V)}^{6/7}\,\le\,D_0\,,
\end{eqnarray*}
where $\,D_0\,$ is a positive constant depending only on the data of the problem.
Moreover, we have
 \begin{eqnarray*}
&&\|\chi_k\|_{L^{7/4}(0,T;L^{7/2}(\Omega))}\,=\,\Bigl[\int_0^T\Bigl(\int_\Omega|\chi_k(x,t)|^{7/2}\,dx\Bigr)^{1/2}dt\Bigr]^{4/7}\\[1mm]
&&=\,\Bigl[\int_0^T\Bigl(\int_\Omega|\chi_k(x,t)|^4\,dx\Bigr)^{1/2}dt\Bigr]^{\frac 1 2\cdot \frac 8 7} \,=\,\|\chi_k\|_{L^2(0,T;L^4(\Omega))}^{8/7}\,.
\end{eqnarray*}

Hence, we can infer from (\ref{eq:3.2}) that for every $k\ge\Phi_0$ it holds the inequality 
\begin{equation}
\label{eq:3.3}
|\!|\!|(\mu-k)^+|\!|\!|\,\le\,k\,D_1\,\|\chi_k\|_{L^2(0,T;L^4(\Omega))}^{8/7}\,,
\end{equation}
where $D_1=D_0/\min\,\{\ve,1\}$, and where the norm $\,|\!|\!|\cdot|\!|\!|\,$ is defined by
$$|\!|\!|v|\!|\!|^2\,:=\,\max_{t\in [0,T]}\|v(t)\|^2_H\,+\,\int_0^T\!\!\!\int_\Omega|\nabla v|^2\,dx\,dt
\quad\forall\,v\in C^0([0,T];H)\cap L^2(0,T;V)\,.$$

Moreover, owing to the continuity of the embedding $V\subset L^4(\Omega)$, there is some
constant $D_2>0$, which only  depends on $\Omega$ and on $T$, such that  
\begin{equation}
\label{eq:3.4}
   \|v\|_{L^2 (0,T; L^4 (\Omega))}\, \leq \,D_2 \,|\!|\!|v|\!|\!|
  \quad \forall\,v\in C^0 ([0,T];H)\cap L^2(0,T;V)\,.
\end{equation}

At this point, we select a strictly increasing sequence $\{k_j\}$
depending on a real parameter $m>1$ as follows:
\begin{equation}
\label{eq:3.5}
  k_j := M \bigl( 2 - 2^{-j} \bigr)
  \quad \mbox{for $j=0,1,\dots$,}\quad\textrm{with}\quad   M := m\, \Phi_0\, .
  \end{equation}
  
Note that $k_0=M>\Phi_0$ and $\lim_{j\to\infty}k_j=2M$.
Then, owing to (\ref{eq:3.3}) 
and (\ref{eq:3.4}), it is not difficult to check that
\begin{eqnarray}
  \label{eq:3.6}
  &&\bigl( k_{j+1} - k_j \bigr) \, \|{\chi_{k_{j+1}}}\|_{L^2 (0,T;
    L^4 (\Omega))}
  \,\leq\, \|(\mu-k_j)^+\|_{L^2 (0,T; L^4 (\Omega))}\nonumber\\[1mm]
  &&\,\leq\, D_2 |\!|\!|(\mu-k_j)^+|\!|\!|
  \,  \leq \,k_j \, D_1\, D_2\, \|\chi_{k_j}\|_{L^2 (0,T; L^4 (\Omega))}^{8/7} .
\end{eqnarray}

Therefore, if we set
$$  S_j := \|\chi_{k_j}\|_{L^2 (0,T; L^4 (\Omega))}
  \quad \mbox{for $j=0,1,\dots$,}$$

then we have
$$  S_{j+1} 
  \,\leq\, \frac {k_j}{k_{j+1}-k_j} \, D_1\, D_2\, S_j^{8/7}
  \,\leq \,4\, D_1\, D_2 \, 2^j\, S_j^{8/7}
  \quad \mbox{for $j=0,1,\dots$}.
  $$
Using [12, Lemma 5.6, p.~95],
we can conclude that $S_j\to0$ as $j\to\infty$, provided that
\begin{equation}
  \label{eq:3.7}
  S_0 \,=\, \|\chi_{k_0}\|_{L^2 (0,T; L^4 (\Omega))}
  \,\leq \,(4\,D_1\,D_2)^{-7}\, 2^{-49}.
\end{equation}

Now recall that $\chi_{k_0}=\chi_M$ and, owing to (\ref{eq:3.5}),
$M>\Phi_0$ and $m=M/\Phi_0$. Also,
$$\chi_M=1<\frac{\mu-\Phi_0}{M-\Phi_0} \quad\mbox{if }\,\mu>M,
\quad\mbox{ and }\,\chi_M=0\quad\mbox{otherwise.}$$
Therefore, using~(\ref{eq:3.3}) and (\ref{eq:3.4}) 
with $k=k_0=M$, we find that
\begin{eqnarray*}
  && S_0 
  \leq \,\frac{1}{M-\Phi_0} \, \|(\mu-\Phi_0)^+\|_{L^2 (0,T; L^4 (\Omega))}
 \, \leq \,\frac{D_2}{M-\Phi_0} \, |\!|\!|(\mu-\Phi_0)^+|\!|\!|
  \\
  && \le \,\frac {D_1\,D_2} {m-1} \, \|\chi_{\Phi_0}\|_{L^2 (0,T; L^4 (\Omega))}^{8/7}
  \leq \frac {D_1\, D_2} {m-1} \, |\Omega|^{\frac14\cdot\frac87} \, 
T^{\frac12\cdot\frac87}.
  \end{eqnarray*}
We are now in a position to choose
$m:=1+D_1\,D_2|\Omega|^{2/7}T^{4/7}(4\,D_1\,D_2)^7\,2^{49}$.
Then, $m>1$ and (\ref{eq:3.7}) is satisfied.
Consequently,
$$  \|\chi_{2M}\|_{L^2 (0,T; L^4 (\Omega))} = \lim_{j\to\infty} S_j = 0,$$
due to Beppo Levi's Monotone Convergence Theorem.
This implies  that $\mu\leq2\,M$ a.e. in $Q\,$,
 and the boundedness of $\mu$ is proved. \qed

\vspace{7mm}
Now that the boundedness condition (\ref{eq:3.1}) is shown, we can prove the following uniqueness and stability result, which corresponds to Theorem~2.2 in \cite{CGPS3}.

\vspace{5mm}
{\bf Theorem 3.2}

(i) {\em\,\,Suppose that} (A1)--(A6) {\em are fulfilled. Then the system}  (\ref{eq:1.1})--(\ref{eq:1.4}) {\em has a unique solution  $(\rho,\mu)$ satisfying} (\ref{eq:2.3})--(\ref{eq:2.6}). 

(ii) {\em Suppose that} (A1){\em ,} (A2){\em ,} (A4) {\em and} (A5) {\em are fulfilled and that the functions 
$\,u_1$, $u_2$ satisfy the conditions} (A3) {\em and} (A6). {\em Moreover, let $(\rho_i,\mu_i)$ be the 
solutions to} (\ref{eq:1.1})--(\ref{eq:1.4}) {\em corresponding to $u_i$, $i=1,2$, and $u:=u_1-u_2$, 
$\rho:=\rho_1-\rho_2$ and $\mu:=\mu_1-\mu_2$. Then we have, for every $t\in [0,T]$,}
\begin{eqnarray}
\label{eq:3.8}
\hskip-.5cm & \displaystyle\max_{0\le s\le t}\,\bigl(\|\mu(s)\|^2_H\,+\|\rho(s)\|_V^2\bigr)\,+\,\int_0^t\!\!\int_\Omega
\bigl(\|\mu(s)\|_V^2\,+\, \|\rho_t(s)\|_H^2\,+\,\|\rho(s)\|_W^2\bigr)\,ds & \nonumber\\[1mm]
\hskip-.5cm&\displaystyle\,\le\,K_1^*\int_0^t \|u(s)\|_{L^2(\Gamma)}^2\,ds\,,&
\end{eqnarray}
{\em with a constant $K_1^*>0$ that only depends on the data of the system.}

\vspace{3mm}
{\em Proof.} \,\quad Obviously, the assertion (i) follows directly from (ii). So we only need to show 
(ii). To this end, observe that by Theorem~3.1 there are constants $M>0$ and $0<r_*<r^*<1$
such that $\,0\le\mu_i\le M\,$ and $\,r_*\le\rho_i\le r^*\,$ a.\,e. in $Q$, for $i=1,2$.
Moreover, the function $\,r\mapsto r-f'(r),$ $r_*\le r\le r^*$, has a Lipschitz constant
$L>0$. Next, we observe that the pair $(\rho,\mu)$ is a solution to the system 
\begin{eqnarray}
\label{eq:3.9}
(\ve+2\rho_1)\mu_t \,+\,2\,\rho\,\mu_{2,t}\,+\,\mu\,\rho_{1,t}+\mu_2\,\rho_t-\Delta
\mu=0\,\quad\mbox{a.\,e. in }\,Q,\\[1mm]
\label{eq:3.10}
\delta\,\rho_t-\Delta\rho=\mu\,-\,(f'(\rho_1)-f'(\rho_2))\,\quad\mbox{a.\,e. in }\,Q,\\[1mm]
\label{eq:3.11}
\frac{\partial \rho}{\partial \bf n}=0\,,\quad\frac{\partial \mu}{\partial \bf n}=
\alpha(u-\mu)\,\quad\mbox{a.\,e. on }\Sigma,\\[1mm]
\label{eq:3.12}
\rho(x,0)= \mu(x,0)=0\,,\quad\mbox{for a.\,e. } \, x\in   \Omega.
\end{eqnarray}

Now observe that $\,2\,\rho_1\,\mu\,\mu_t=\bigl(\rho_1\,\mu^2\bigr)_t\,-\,\mu^2\,\rho_{1,t}$. Hence, if we test 
(\ref{eq:3.9}) by $\mu$ then we obtain, using Young's inequality, that for every $t\in [0,T]$ it holds
\begin{eqnarray}
\label{eq:3.13}
&&\int_\Omega \left(\frac \ve 2\,+\,\rho_1(t)\right)\mu^2(t)\,dx\,+\,\int_0^t\!\!\int_\Omega |\nabla\mu|^2\,
dx\,ds\,+\,\int_0^t\!\!\int_\Gamma \alpha\,|\mu|^2\,d\sigma\,ds\nonumber\\[1mm]
&&\le\,C\int_0^t\!\!\int_\Gamma |u|^2\,d\sigma\,ds\,+\,
\int_0^t\!\!\int_\Gamma\frac\alpha 2\,|\mu|^2\,d\sigma\,ds
\,\nonumber\\[1mm]
&&\quad+\,
\int_0^t\!\!\int_\Omega |\mu|\,\bigl(2\,|\rho|\,|\mu_{2,t}|\,+\,|\mu_2|\,|\rho_t|\bigr)
dx\,ds\,.
\end{eqnarray}

We have, owing to the continuity of the embedding $\heins\subset L^4(\Omega)$ and to Young's inequality,
\begin{eqnarray}
\label{eq:3.14}
&&\int_0^t\!\!\int_\Omega 2\,|\mu|\,|\rho|\,|\mu_{2,t}|\,dx\,ds\,\le\,
C\int_0^t\|\mu_{2,t}(s)\|_H^2\,\|\mu(s)\|_{L^4(\Omega)}\|\rho(s)\|_{L^4(\Omega)}\,ds\nonumber\\[1mm]
&&\le\,\gamma\int_0^t \|\mu(s)\|_V^2\,ds\,+\,\frac C {\gamma} \int_0^t\|\mu_{2,t}(s)\|_H^2\,\|\rho(s)\|_V^2\,ds\,,
\end{eqnarray}
where, owing to (\ref{eq:2.4}), the mapping $s\mapsto \|\mu_{2,t}(s)\|_H^2$ belongs to $L^1(0,T)$. Moreover, we also
have $\mu_2\in\liq$, and thus
\begin{eqnarray}
\label{eq:3.15}
&&\int_0^t\!\!\int_\Omega |\mu|\,|\mu_2|\,|\rho_t|\,dx\,ds\,\le\,
C\int_0^t\|\rho_t(s)\|_H\,\|\mu(s)\|_H\,ds\nonumber\\[1mm]
&&\le\,\gamma\int_0^t \|\rho_t(s)\|_H^2\,ds\,+\,\frac C {\gamma} \int_0^t\|\mu(s)\|_H^2\,ds\,.
\end{eqnarray}

Next, we add $\rho$ on both sides of Eq. (\ref{eq:3.10}) and test the resulting equation by $\rho_t$. 
Invoking Young's inequality, it is easily seen that, for every $t\in [0,T]$,
\begin{eqnarray}
\label{eq:3.16}
&&\delta\int_0^t\|\rho_t(s)\|_H^2\,ds\,+\,\|\rho(t)\|_V^2\nonumber\\[1mm]
&&\le\,\gamma\int_0^t\|\rho_t(s)\|_H^2\,ds\,+\,\frac C {\gamma} \int_0^t\bigl(\|\mu(s)\|_H^2\,+\,
L^2\,\|\rho(s)\|_H^2\bigr)\,ds\,.
\end{eqnarray}

Now we can combine (\ref{eq:3.13})--(\ref{eq:3.16}). Choosing $\gamma>0$ sufficiently small, and 
applying Gronwall's lemma, we see that (\ref{eq:3.8}) is satisfied. \qed  

\vspace{5mm}
The stability estimate (\ref{eq:3.8}) can be improved if further regularity is assumed for $f$. 
The following result is a counterpart of Lemma 3.1 in \cite{CGPS5}. We remark at this place that (\ref{eq:2.3}) implies,
in particular, that $\rho$ is weakly continuous 
as a mapping from $[0,T]$ into~$W$, which justifies the formulation of the estimate (3.17) below.

\vspace{5mm}
{\bf Theorem 3.3} \quad\,\,{\em Suppose that the assumptions of Theorem~3.2,(ii) are satisfied, and assume that} 

\vspace{1mm}
(A7) \quad $f\in C^3(0,1)$.

\vspace{1mm} 
{\em Then we have, for every $t\in [0,T]$,}
\begin{eqnarray}
\label{eq:3.17}
&&\max_{0\le s\le t}\,\bigl(\|\mu(s)\|^2_V\,+\|\rho_t(s)\|_V^2\,+\,\|\rho(s)\|_W^2\bigr)
\nonumber\\[1mm]
&& \quad {}+\int_0^t
\bigl(\|\mu_t(s)\|_H^2\,+\,
\|\rho_t(s)\|_W^2\bigr)\,ds
\nonumber\\[1mm]
&&\,\le\,K_2^*
\left\{
\|u(0)\|_{L^2(\Gamma)}^2 
+\int_0^t \left(\|u(s)\|_{L^2(\Gamma)}^2\,+\,\|u_t(s)\|_{L^2(\Gamma)}^2\right)\,ds
\right\}
\qquad
\end{eqnarray}
{\em with a constant $K_2^*>0$ that only depends on the data of the system.}

\vspace{3mm}
{\bf Remark 3.4} \quad
We note that 
$\|u(0)\|_{L^2(\Gamma)}^2\leq I(t)/\mathop{\rm max}\{1,t\}$
where $I(t)$ denotes the last integral of (\ref{eq:3.17}).
It follows that $\|u(0)\|_{L^2(\Gamma)}^2$ can be dropped
if one pretends (\ref{eq:3.17}) just for $t=T$.

\vspace{3mm}
{\em Proof of Theorem~3.3}. \,We closely follow the lines of the proof of Lemma 3.1 in \cite{CGPS5}. Since the proof
given there carries over to our situation with minor changes, we can afford to be brief. 
First, observe that by Theorem~3.1 there are constants $M>0$ and $0<r_*<r^*<1$
such that $\,0\le\mu_i\le M\,$ and $\,r_*\le\rho_i\le r^*\,$ a.\,e. in $Q$, for $i=1,2$.
Next, we recall that the pair $(\rho,\mu)$ is a solution to the system 
(\ref{eq:3.9})--(\ref{eq:3.12}). 
We test Eq. (\ref{eq:3.9}) formally by $\mu_t$. 
It then follows, with the use of Young's inequality, that
\begin{eqnarray}
\label{eq:3.18}
&& \ve\int_0^t\|\mu_t(s)\|^2_H\,ds \,+\,\frac 1 2 \|\nabla\mu(t)\|_H^2\,
+\,\frac 1 2\int_\Gamma \alpha\,|\mu(t)|^2\,d\sigma\nonumber\\[1mm]
&&\le\,\int_0^t\!\!\int_\Gamma\alpha\,u\,\mu_t\,d\sigma\,ds\,\nonumber\\[1mm]
&& \quad +\! \int_0^t\!\!\int_\Omega\left(2|\rho|\,|\mu_{2,t}|\,+\,|\mu|\,|\rho_{1,t}|\,+\,|\mu_2|
\,|\rho_t|\right)|\mu_t|\,dx\,ds\,.
\end{eqnarray}
Now, by virtue of integration by parts with respect to $t$, and invoking (\ref{eq:3.8}), Young's inequality and the 
trace theorem,
\begin{eqnarray}
\label{eq:3.19}
\hspace*{-15mm}&&\Bigl|\int_0^t\!\!\int_\Gamma\alpha\,u\,\mu_t\,d\sigma\,ds\Bigr|
\nonumber\\[1mm]
\hspace*{-15mm}&&\leq\gamma\int_\Gamma\alpha|\mu(t)|^2\,d\sigma\,
+\,\frac C\gamma \,\int_\Gamma |u(t)|^2\,d\sigma\,
+ \,\frac C\gamma \,\int_0^t\!\!\int_\Gamma|u_t|\,|\mu|\,d\sigma\,ds
\nonumber\\[1mm]
\hspace*{-15mm}&&\leq C\,\gamma\,\|\mu(t)\|_V^2\,
\,+\,\frac C\gamma \int_\Gamma |u(0)|^2\,d\sigma
\,+\,\frac C\gamma\int_0^t\!\!\int_\Gamma( |u|^2\,+\,|u_t|^2)\,d\sigma\,ds\,.
\end{eqnarray}
 
Employing almost exactly the same arguments as in the proof of Lemma 3.1 in \cite{CGPS5}
(the minor necessary changes are left as an easy exercise to the reader), and taking advantage of (\ref{eq:3.8}),  we conclude the estimate (where $\gamma>0$ is arbitrary)
\begin{eqnarray}
\label{eq:3.20}
&&\int_0^t\!\!\int_\Omega\left(2|\rho|\,|\mu_{2,t}|\,+\,|\mu|\,|\rho_{1,t}|\,+\,|\mu_2|
\,|\rho_t|\right)|\mu_t|\,dx\,ds\nonumber\\[1mm]
&&\le\,3\,\gamma\int_0^t\|\mu_t(s)\|_H^2\,ds\,+\,\frac C \gamma \int_0^t\|\rho_{1,t}(s)\|_V^2\,\|\mu(s)\|_V^2
\,ds\nonumber\\[1mm]
&&\quad +\,\frac C \gamma \int_0^t\|\mu_{2,t}(s)\|_H^2\,\|\Delta\rho(s)\|_H^2 \,ds\,
\,+\,\frac C\gamma
\int_0^t\!\!\int_\Gamma |u|^2\,d\sigma\,ds\,.
\qquad\quad
\end{eqnarray}

Next, we test (\ref{eq:3.10}) formally by $\,-\Delta\rho_t$. By the same token as in the proof of
Lemma 3.1 in \cite{CGPS5}, we deduce for arbitrary $\gamma>0$ the estimate
\begin{eqnarray}
\label{eq:3.21}
&&\delta\int_0^t\|\nabla\rho_t(s)\|_H^2\,ds\,+\,\frac 1 4 \,\|\Delta\rho(t)\|_H^2
\,\le\,\gamma\int_0^t\|\mu_t(s)\|^2_H\,ds\nonumber\\[1mm]
&&+\,\frac C{\gamma}\int_0^t\left(1\,+\,\|\rho_{2,t}(s)\|_V^2\right)\,\|\Delta\rho(s)\|_H^2\,ds
\,+\,C\int_0^t\!\!\int_\Gamma |u|^2\,d\sigma\,ds\,.\qquad
\end{eqnarray}

Now observe that, owing to Theorem~2.1, the mappings $\,s\mapsto \|\rho_{i,t}(s)\|_V^2\,$, $i=1,2$, and 
$\,s\mapsto\|\mu_{2,t}(s)\|^2_H\,$ all belong to $L^1(0,T)$.
Hence, combining the estimates (\ref{eq:3.18})--(\ref{eq:3.21}), adjusting $\gamma>0$ sufficiently small, and invoking Gronwall's lemma, we can conclude that for every $t\in [0,T]$ it holds
\begin{eqnarray}
\label{eq:3.22}
&&\int_0^t\left(\|\nabla\rho_t(s)\|_H^2\,+\,\|\mu_t(s)\|_H^2\right)ds \,+\,\max_{0\le s\le t}\,
\left(\|\mu(s)\|_V^2\,+\,\|\rho(s)\|_W^2\right)\nonumber\\[1mm]
&&\le\,C
\,\left\{
\|u(0)\|_{L^2(\Gamma)}^2 
+ \int_0^t\!\!\int_\Gamma \left(|u|^2\,+\,|u_t|^2\right)\,d\sigma\,ds
\right\}
\,.\qquad
\end{eqnarray}

Next, we formally differentiate (\ref{eq:3.10}) with respect to $t$, and obtain
\begin{equation}
\label{eq:3.23}
\delta\rho_{tt}-\Delta\rho_t=\mu_t\,-\,f''(\rho_1)\,\rho_t -(f''(\rho_1)-f''(\rho_2))\,\rho_{2,t}\,,
\end{equation}
with zero initial and Neumann boundary conditions for $\rho_t$. Hence, testing 
(\ref{eq:3.23}) by $\rho_t$, invoking
Young's inequality, and recalling (\ref{eq:3.8}) and (\ref{eq:3.22}), we find that 
\begin{eqnarray}
\label{eq:3.24}
&&\frac {\delta} {2}\, \|\rho_t(t)\|_H^2\,+\,\int_0^t\!\!\|\nabla\rho_t(s)\|_H^2\,ds
\nonumber\\
&&\le C\,
\left\{
\|u(0)\|_{L^2(\Gamma)}^2 
+ \int_0^t\!\!\int_\Gamma \left(|u|^2\,+\,|u_t|^2\right)\,d\sigma\,ds
\right\}
\nonumber\\
&& \quad {}+\int_0^t\!\!\int_\Omega|\rho_{2,t}|\,|f''(\rho_1)-f''(\rho_2)|\,|\rho_t|\,dx\,ds\,. 
\end{eqnarray}

Moreover, using H\"older's and Young's inequalities, (A7) and (\ref{eq:3.8}), 
we see that
\begin{eqnarray}
\label{eq:3.25}
&&\int_0^t\!\!\int_\Omega|\rho_{2,t}|\,|f''(\rho_1)-f''(\rho_2)|\,|\rho_t|\,dx\,ds\nonumber\\
&&\,\le\,C\int_0^t\|\rho_{2,t}(s)\|_{L^4(\Omega)}\,\|\rho(s)\|_{L^4(\Omega)}\|\rho_t(s)\|_H\,ds\nonumber\\
&&\,\le\,C\,\Bigl(\int_0^t\|\rho_t(s)\|_H^2\,ds\,+\,\max_{0\le s\le t}\|\rho(s)\|_V^2
\int_0^t\|\rho_{2,t}(s)\|_V^2\,ds\Bigr)\nonumber\\
&&\,\le\,C\int_0^t\!\!\int_\Gamma(|u|^2\,+\,|u_t|^2)\,d\sigma\,ds\,.
\end{eqnarray}

\vspace{2mm}
Finally, we test (\ref{eq:3.23}) by $\,-\Delta\rho_t$. Using Young's inequality and
(\ref{eq:3.22}), we find that

\begin{eqnarray}
\label{eq:3.26}
&&\frac {\delta} 2 \,\|\nabla\rho_t(t)\|_H^2\,+\, \int_0^t\|\Delta\rho_t(s)\|_H^2\,ds
\nonumber\\[1mm]
&&\,\le\,\gamma\int_0^t\|\Delta\rho_t(s)\|_H^2\,ds
\,+\,\frac{C}{\gamma}
\left\{
\|u(0)\|_{L^2(\Gamma)}^2 
+ \int_0^t\!\!\int_\Gamma \left(|u|^2\,+\,|u_t|^2\right)\,d\sigma\,ds
\right\}
\nonumber\\[1mm]
&& \quad + \int_0^t\!\!\int_\Omega|\rho_{2,t}|\,|f''(\rho_1)-f''(\rho_2)|\,|\Delta\rho_t|\,dx\,ds
\nonumber\\[1mm]
&&\le \,2\gamma\!\int_0^t\!\!\|\Delta\rho_t(s)\|_H^2\,ds
\,+ \,\frac{C}{\gamma}
\left\{
\|u(0)\|_{L^2(\Gamma)}^2 
+ \int_0^t\!\!\int_\Gamma \left(|u|^2\,+\,|u_t|^2\right)\,d\sigma\,ds
\right\}
\nonumber\\[1mm]
&&\quad {} + \frac C{\gamma}
\,\max_{0\le s\le t}\|\rho(s)\|_V^2
\int_0^t\|\rho_{2,t}(s)\|_V^2\,ds
\nonumber\\[1mm]
&&\le \,2\gamma\!\int_0^t\!\!\|\Delta\rho_t(s)\|_H^2\,ds
\nonumber\\[1mm]
&& \quad {} +\, \frac C{\gamma}
\left\{
\|u(0)\|_{L^2(\Gamma)}^2 
+ \int_0^t\!\!\int_\Gamma \left(|u|^2\,+\,|u_t|^2\right)\,d\sigma\,ds
\right\}.
\end{eqnarray}

Choosing $\gamma>0$ appropriately small, we can infer that the estimate (\ref{eq:3.17}) is 
in fact true. This concludes the proof.
\qed

\section{An optimal boundary control problem}
\setcounter{equation}{0}
In this section, we consider the following optimal boundary control problem:

\vspace{3mm}
{\bf (CP)} \quad Minimize the cost functional
\begin{eqnarray}
\label{eq:4.1}
J(u,\rho,\mu)&\!\!=\!\!&\frac 12 \int_\Omega|\rho(x,T)-\rho_T(x)|^2\,dx\,+\,\frac{\beta_1} 2 \int_0^T\!\!
\int_\Gamma|u(x,t)|^2\,d\sigma\,dt\nonumber\\[1mm]
&&+\,\frac{\beta_2} 2 \int_0^T\!\!\int_\Omega|\mu(x,t)-\mu_T(x,t)|^2\,dx\,dt
\end{eqnarray}

subject to the state system (\ref{eq:1.1})--(\ref{eq:1.4}) and to the control constraints
\begin{eqnarray}
\label{eq:4.2}
&&u\in\uad := \left\{v\in H^1(0,T;L^2(\Gamma))\cap \lis:\,U_1\le v\le U_2 \,\,\,\,\mbox{a.\,e. on }\,
\Sigma, \right.\nonumber\\[2mm]
&&\hskip2cm \,\,\,\left.\|v_t\|_{L^2(0,T;L^2(\Gamma))}\le R\right\}\,.
\end{eqnarray}

In this connection, we require that the hypotheses (A1)--(A7) be
satisfied. In addition, 
we postulate:

\vspace{3mm}
(A8) \quad $R>0$, $\beta_i\ge 0$, $i=1,2$, $\rho_T\in\lzo$,  $\mu_T\in
\lzq$, $U_1,U_2\in$ 
\linebreak
\hspace*{11mm} 
$\lis$, and 
there are constants
$0< u_* <u^*<+\infty$ such that
\begin{equation}
\label{eq:4.3}
u_* \le U_1\le U_2\le  u^*\quad\mbox{a.\,e. on }\,\Sigma\,.
\end{equation}

\vspace{5mm}
In what follows, we denote 
\begin{equation*}
\mathcal{X}:=\heg\cap\lis\,, \quad \|u\|_{\mathcal{X}}:=\|u\|_{\heg}\,+\,\|u\|_{\lis} \,,
\end{equation*}
where $\|\cdot\|_{\heg}$ denotes the standard norm in $\heg$.
Obviously, $\uad$ is a nonempty, bounded, closed and convex subset 
of $\mathcal{X}$, and $\uad$ is contained in the open set $\mathcal{U}\subset\mathcal{X}$ given by
\begin{equation*}
\mathcal{U}\!:=\!\left\{v\in \mathcal{X}: \, \ \frac 1 2\, u_* < \mbox{ess\,inf}\,v \, , \ \  \mbox{ess\,sup}\, v<\frac 3 2\, u^*
\,, \,\,\ \|v_t\|_{L^2(0,T;L^2(\Gamma))}<R\!+\!1 \right\}.
\end{equation*}
By Theorem~3.3, the control-to-state mapping $\,u\mapsto S(u):=(\rho,\mu)$ is 
Lip\-schitz continuous as a mapping from the set
$\mathcal{U}\subset \mathcal{X}$ into the space
$$\left(H^1(0,T;W)\cap C^1([0,T];V)\right)\times \left(H^1(0,T;H)\cap C^0([0,T];V)\right).$$ 
We may without loss of generality assume (by possibly taking a larger $K_2^*$) that (\ref{eq:3.17}) is valid on the whole set
$\mathcal{U}$ with the same constant $K^*_2>0$. It also follows from Theorem~3.1 that there exist constants $\mu^*>0$ and
$0<r_*<r^*<1$ such that for every $u\in \mathcal{U}$ it holds 
\begin{equation}
\label{eq:4.4}
0\le\mu\le \mu^*\quad\mbox{and }\, 0<r_*\le \rho\le r^*<1 \quad\mbox{a.\,e. in }\,Q,
\end{equation}
where $(\rho,\mu)=S(u)$. Moreover, a closer inspection of the proof of Theorem~2.1 reveals that there is a constant
$K_3^*>0$ such that we have, for any $u\in \mathcal{U}$,
\begin{eqnarray}
\label{eq:4.5}
&&\|\rho\|_{W^{1,\infty}(0,T;H)\cap H^1(0,T;V)\cap L^\infty(0,T;W)}\nonumber\\[2mm]
&&+\,\|\mu\|_{H^1(0,T;H)\cap C^0([0,T];V)\cap L^2(0,T;\hth)\cap \liq}\,\le\,K^*_3\,.
\end{eqnarray}

\vspace{3mm}
{\bf Remark 4.1} \quad Thanks to (\ref{eq:4.4}) and to $f\in C^3(0,1)$, it holds $f'(\rho)\in\liq$. Also, by
the embedding $V\subset L^6(\Omega)$, we have $\mu\in C^0([0,T];L^6(\Omega))$. Notice also that (\ref{eq:2.3}) implies, in particular, that $\rho$ is continuous from $[0,T]$ to $H^s(\Omega)$ for all $s<2$; thus, since $H^s(\Omega)\subset
C(\overline\Omega)$ for $s>3/2$, 
we also have $\rho\in C(\overline Q)$. 
Therefore, possibly choosing a larger constant $K_3^*$, we may without loss of generality assume that
\begin{equation}
\label{eq:4.6}
\|\rho\|_{C(\overline Q)}\,+\,\|\mu\|_{C^0([0,T];L^6(\Omega))}\,+\,\|\rho_t\|_{L^2(0,T;L^6(\Omega))}\,\le
\,K^*_3\,\quad \forall\,u\in \mathcal{U}.
\end{equation}

\vspace{3mm}
{\bf Remark 4.2} \quad  The mathematical literature on control problems for phase field systems is scarce and
usually restricted to the so-called {\em Caginalp model} of phase transitions (see, e.\,g., 
\cite{HJ}, \cite{Hei}, \cite{HeiTr}, \cite{Tr}, and the references given there). More general,
thermodynamically consistent phase field models were the subject of \cite{LS07}. In \cite{CGPS5}, the 
present authors investigated a control problem for the system (\ref{eq:1.1})--(\ref{eq:1.4}) 
with distributed controls. Since many  of the arguments employed in \cite{CGPS5} carry over to the boundary control
considered here, we can afford to be sketchy in some of the proofs in the following exposition.

\subsection{Existence}
We begin our discussion of the control problem {\bf (CP)} with the following existence result:

\vspace{5mm}
{\bf Theorem 4.3} \quad {\em Suppose that the conditions} (A1)--(A8) {\em are satisfied. Then
the optimal control problem} {\bf (CP)} {\em has a solution} $\overline u \in \uad$.

\vspace{3mm}
{\em Proof.} Let $\{u_n\}\subset \uad$ be a minimizing sequence for {\bf (CP)}, and let $\{(\rho_n,
\mu_n)\}$ be the sequence of the associated solutions to (\ref{eq:1.1})--(\ref{eq:1.4}). We then can infer from 
(\ref{eq:4.5}) the 
existence of a triple $(\bar u,\bar\rho,\bar\mu)$ such that, for a suitable subsequence again indexed by
$n$, we have 
\begin{eqnarray*}
\hskip-.5cm&&u_n\to \bar u\,\quad\mbox{weakly star in }\,\heg\cap L^\infty(\Sigma),\\[1mm]
\hskip-.5cm&&\rho_n\to\bar\rho\,\quad\mbox{weakly star in }\,W^{1,\infty}(0,T;H)\cap H^1(0,T;V)\cap L^\infty(0,T;W),\\[1mm]
\hskip-.5cm&&\mu_n\to\bar\mu \,\quad\mbox{weakly star in }\,H^1(0,T;H)\cap L^\infty([0,T];V)\cap L^2(0,T;\hth).
\end{eqnarray*}
Clearly, we have that $\bar u\in \uad$. Moreover, by virtue of the Aubin-Lions lemma 
(cf. \cite[Thm. 5.1, p. 58]{Aubin}) and similar compactness results
(cf. \cite[Sect. 8, Cor. 4]{Simon}), we also have the strong convergences
\begin{eqnarray*}
&&\rho_n\to \bar\rho\,\quad\mbox{strongly in }\,C^0([0,T];H^s(\Omega)) \,\quad\mbox{for
all }\,s<2,\\[1mm]
&&\mu_n\to \bar\mu\,\quad\mbox{strongly in }\,C^0([0,T];H)\cap L^2(0,T;V).
\end{eqnarray*}
From this we infer, possibly selecting another subsequence again indexed by $n$, that 
$\rho_n\to\bar\rho $ pointwise a.\,e. (actually, uniformly) in $Q$. In particular,
$r_*\le \bar\rho\le r^*$ a.\,e. in $Q$ and, since $f\in C^2(0,1)$, also
$f'(\rho_n)\to f'(\bar\rho)$  strongly in $L^2(Q)$. 
Now notice that the above convergences imply, in particular, that
\begin{eqnarray*}
&&\rho_n\to \bar\rho\quad \mbox{strongly in}\,C^0([0,T];L^6(\Omega)),\\[1mm] 
&&\partial_t\rho_n\to \partial_t\bar\rho \quad\mbox{weakly in }\,L^2(0,T;L^4(\Omega)),\\[1mm]
&&\mu_n\to \bar\mu\quad\mbox{strongly in }\,
L^2(0,T;L^4(\Omega)),\\[1mm]
&&\partial_t\mu_n\to \partial_t\bar\mu\quad\mbox{weakly in }\,L^2(Q).
\end{eqnarray*}
From this, it is easily verified that
\begin{eqnarray*}
&&\mu_n\,\partial_t\rho_n\to \bar\mu\,\partial_t\bar\rho\,\quad\mbox{weakly in }\,
L^1(0,T;H), \\[1mm]
&&\rho_n\,\partial_t\mu_n\to\bar\rho\,\partial_t\bar\mu\,\quad\mbox{weakly in }\,
L^2(0,T;L^{3/2}(\Omega)).
\end{eqnarray*}
In summary, if we pass to the limit as $n\to\infty$ in the state equations 
(\ref{eq:1.1})--(\ref{eq:1.4})
written for the triple $(u_n,\rho_n,\mu_n)$, we find that $(\bar\rho,\bar\mu)=S(\bar u)$,
that is, the triple $(\bar u,\bar\rho,\bar\mu)$ is admissible for the control problem
{\bf (CP)}. From the weak sequential lower semicontinuity of the cost functional $J$
it finally follows that $\bar u$, together with $(\bar\rho,\bar\mu)=S(\bar u)$, is
a solution to {\bf (CP)}. This concludes the proof. \qed

\vspace{4mm}
{\bf Remark 4.4} \quad It can be shown that this existence result holds for much more general cost
functionals. All we need is that $J$ enjoy appropriate weak sequential
lower semicontinuity properties that match the above weak convergences.

\vspace{4mm}
{\bf Remark 4.5} \quad  Since the state component $\rho$ is continuous on $\overline Q$, the existence 
result remains valid if suitable pointwise state constraints for $\rho$
are added (provided the admissible set is not empty). 

\subsection{Necessary optimality conditions}

In this section, we derive the first-order necessary conditions of optimality for
problem {\bf (CP)}. To this end, we first show that the control-to-state operator $S:u\mapsto (\rho,\mu)$ is
Fr\'{e}chet differentiable as a mapping from $\mathcal{U}\subset \mathcal{X}$ into the Banach space  $\left(\mathcal{Y},
\|\,\cdot\,\|_\mathcal{Y}\right)$, where
\begin{eqnarray*}
\mathcal{Y}&\!\!:=\!\!&\left(H^1(0,T;H)\cap C^0([0,T];V)\cap L^2(0,T;W)\right)\nonumber\\[1mm]
&&\times \left(C^0([0,T];H)\cap L^2(0,T;V)\right) .
\end{eqnarray*}
          
\subsubsection{The linearized system}

Let $\bar u\in \mathcal{U}$ and $h\in \mathcal{X}$ be given and $(\bar\rho,\bar\mu)=S(\bar u)$. As a preparatory step, we consider the following system, which is obtained by linearizing the system 
(\ref{eq:1.1})--(\ref{eq:1.4}) at $(\bar\rho,\bar\mu)$:
\begin{eqnarray}
\label{eq:4.7}
(\ve+2\bar\rho)\,\eta_t-\Delta\eta+2\,\bar\mu_t\,\xi+\bar\mu\,\xi_t+\bar\rho_t\,\eta=0\,\quad
\mbox{a.\,e. in }\,Q,\\[1mm]
\label{eq:4.8}
\delta\,\xi_t-\Delta\xi=-f''(\bar\rho)\,\xi\,+\,\eta\,\quad
\mbox{a.\,e. in }\,Q,\\[1mm]
\label{eq:4.9}
\,\,\frac{\partial\xi}{\partial \bf n}=0\,,\quad\frac{\partial\eta}{\partial \bf n}=\alpha
(h-\eta)\,\quad
\mbox{a.\,e. on }\Sigma,\\[1mm]
\label{eq:4.10}
\xi(x,0)=\eta(x,0)=0 \,\,\quad \mbox{for a.\,e. } \, x\in \Omega.
\end{eqnarray}

We expect for the Fr\'echet derivative $DS(\bar u)$ at $\bar u$ (if it exists) that $(\xi,\eta)=DS(\bar u)h$, 
provided that (\ref{eq:4.7})--(\ref{eq:4.10}) admits a unique solution
$(\xi,\eta)$.
In view of (\ref{eq:2.3}), (\ref{eq:2.4}), and (\ref{eq:3.1}), 
we can guess the regularity of $\,\xi\,$ and $\,\eta\,$:
\begin{eqnarray}
\label{eq:4.11}
&&\xi\in H^1(0,T;H)\cap C^0([0,T];V)\cap L^2(0,T;W)\cap L^\infty(Q),\\[1mm]
\label{eq:4.12}
&&\eta\in  H^1(0,T;H)\cap C^0([0,T];V)\cap L^2(0,T;\hth).
\end{eqnarray}
Notice that also in this case we cannot expect that $\eta(t)\in H^2(\Omega)$ a.e.\ in $(0,T)$
due to the low space regularity of~$h$,
and we could repeat Remark~2.2 here.
Nevertheless, 
if (\ref{eq:4.11}) and (\ref{eq:4.12}) hold, then the collection of source terms in 
(\ref{eq:4.7}), i.\,e., the part $-2\,\bar\mu_t\,\xi-\bar\mu\,\xi_t-\bar\rho_t\,\eta$,
belongs to $L^2(Q)$,
whereas the regularity (\ref{eq:4.12}) for $\eta$ allows us to conclude from (\ref{eq:4.8})
that also $\xi\in C(\overline Q)$ (by applying maximal parabolic regularity theory, see, e.\,g.,
\cite[Thm. 6.8]{Grie} or \cite[Lemma 7.12]{Tr}).

In fact,  $\xi$ is even more regular: indeed, we may differentiate (\ref{eq:4.8}) with respect to $t$ to find that
\begin{equation} 
\label{eq:4.13}
\delta\xi_{tt}-\Delta\xi_t=-f'''(\bar\rho)\,\bar\rho_t\,\xi-f''(\bar\rho)\,\xi_t+\eta_t\,,
\end{equation}
with zero initial and Neumann boundary conditions for $\,\xi_t$. Since the right-hand side of (\ref{eq:4.13})
belongs to $L^2(Q)$, we may test by any of the functions
$\xi_t$, $\xi_{tt}$, and $\,-\Delta\xi_t$, to obtain that even
\begin{equation}
\label{eq:4.14}
\xi\in H^2(0,T;H)\cap C^1([0,T];V)\cap H^1(0,T;W)\,.
\end{equation} 
Notice, however, that this fact
 has no bearing on the regularity of $\eta$, since the coefficient $\bar\mu_t$ in (\ref{eq:4.7})
only belongs to $L^2(Q)$.  

\vspace{3mm} The following well-posedness result resembles Proposition~3.2 in \cite{CGPS5}.

\vspace{3mm} 
{\bf Proposition 4.6} \quad {\em Suppose that} (A1)--(A8) {\em are fulfilled. Then
 the system} (\ref{eq:4.7})--(\ref{eq:4.10})
 {\em has a unique solution 
$(\xi,\eta)$ satisfying} (\ref{eq:4.12}), (\ref{eq:4.14}), {\em  and
\begin{eqnarray}
\label{eq:4.15}
&\hspace*{-0.8cm}&\|\xi\|_{H^2(0,T;H)\cap C^1([0,T];V)\cap H^1(0,T;W)} 
\,+\,\|\eta\|_{H^1(0,T;H)\cap C^0([0,T];V)\cap L^2(0,T;\hth)}\nonumber\\[1mm] 
&\hspace*{-0.8cm}&\le \,K_4^*\,\|h\|_{H^1(0,T;L^2(\Gamma))}\,,
\end{eqnarray} 
with a constant $K_4^*>0$ that is independent of the choice of $\bar u\in \mathcal{U}$ and
$h\in \mathcal{X}$.}     

\vspace{4mm}
{\bf Remark 4.7} \quad It follows from Proposition 4.6, in particular, that 
the linear mapping $h\mapsto (\xi,\eta)$ is continuous from $\mathcal{X}$ into $\mathcal{Y}$.

\vspace{5mm}
{\em Proof.} \,We follow the lines of the proof of our previous existence results 
and proceed in a series of steps. 

\underline{Step 1: \,Approximation.} \,\quad As in the proof of Theorem~2.1, we 
use an approximation technique based on a delay in the right-hand side of (\ref{eq:4.8}). To this
end, for $\tau>0$ we resume the definition of the translation operator $\mathcal{T}_\tau:\,L^1(0,T;H)$ $\to L^1(0,T;H)$ by putting, for every $v\in L^1(0,T;H)$ and almost every $t\in (0,T)$, 
\begin{equation}
\label{eq:4.16}
(\mathcal{T}_\tau v)(t)= v(t-\tau)\,\mbox{ if }\,t\ge\tau, \quad\,\mbox{and }\,\,\,
(\mathcal{T}_\tau v)(t)= 0\,\mbox{ if }\,t<\tau. 
\end{equation}
Notice that, for any $v\in L^2(Q)$ and any $\tau>0$, we obviously have 
$\,\|\mathcal{T}_\tau v\|_{L^2(Q)}$ $\le \|v\|_{L^2(Q)}$.

Then, for any fixed $\tau>0$, we look for functions $(\xi^\tau, \eta^\tau)$, which satisfy (\ref{eq:4.11}) and (\ref{eq:4.12}) and the system: 
\begin{eqnarray}
\label{eq:4.17}
(\ve+2\bar\rho)\,\eta^\tau_t-\Delta\eta^\tau+2\,\bar\mu_t\,\xi^\tau
+\bar\mu\,\xi^\tau_t+\bar\rho_t\,\eta^\tau=0\,\quad
\mbox{a.\,e. in }\,Q,\\[1mm]
\label{eq:4.18}
\delta\,\xi^\tau_t-\Delta\xi^\tau+f''(\bar\rho)\,\xi^\tau=\mathcal{T}_\tau\eta^\tau\,\quad
\mbox{a.\,e. in }\,Q,\\[1mm]
\label{eq:4.19}
\,\,\frac{\partial\xi^\tau}{\partial \bf n}=0\,,\quad\frac{\partial\eta^\tau}{\partial \bf n}=
\alpha(h-\eta^\tau)\,\quad
\mbox{a.\,e. on }\Sigma,\\[1mm]
\label{eq:4.20}
\xi^\tau(x,0)=\eta^\tau(x,0)=0 \,\,\quad\mbox{for a.\,e. }\,x\in \Omega .
\end{eqnarray}
Precisely, we choose for $\tau>0$ the discrete values $\tau=T/N$, where $N\in\nz$ is arbitrary, and put
$t_n=n\,\tau$, $0\le n\le N$, and  $I_n=(0,t_n)$. For $1\le n\le N$, we solve the problem
\begin{eqnarray}
\label{eq:4.21}
(\ve+2\bar\rho)\,\eta_{n,t}-\Delta\eta_n+2\,\bar\mu_t\,\xi_n+\bar\mu\,\xi_{n,t}+\bar\rho_t\,\eta_n=0\,\quad
\mbox{a.\,e. in }\,\Omega\times I_n,\quad\\[1mm]
\label{eq:4.22}
\frac{\partial\eta_n}{\partial \bf n}=\alpha(h-\eta_n)\,\quad\mbox{a.\,e. on }\,\Gamma\times I_n, \,\quad
\eta_n(x,0)=0 \,\,\quad\mbox{for a.\,e. }\,x\in \Omega,\quad \\[1mm]
\label{eq:4.23}
\delta\,\xi_{n,t}-\Delta\xi_n+f''(\bar\rho)\,\xi_n=\mathcal{T}_\tau\eta_n\,\quad
\mbox{a.\,e. in }\,\Omega\times I_n,\quad\\[1mm]
\label{eq:4.24}
\frac{\partial\xi_n}{\partial \bf n}= 0 \,\quad\mbox{a.\,e. on }\,\Gamma\times I_n,\,\quad
\xi_n(x,0)=0 \,\,\quad \mbox{for a.\,e. }\,x\in \Omega,\quad
\end{eqnarray}
where the variables $\eta_n$ and $\xi_n$, defined on $I_n$, have obvious meaning.
Here, $\mathcal{T}_\tau$ acts on functions that are not defined on the entire interval $(0,T)$; however,
for $n>1$ it is still defined by (\ref{eq:4.16}), while for $n=1$ we simply 
put $\mathcal{T}_\tau\eta_n=0$.
Notice that whenever the pairs $(\xi_k,\eta_k)$ with
\begin{eqnarray}
\label{eq:4.25}
\xi_k\in H^1(I_k;H)\cap C^0(\bar I_k;V)\cap L^2(I_k;W)\cap C(\overline {\Omega\times I_k}),\\
\label{eq:4.26}
\eta_k\in  H^1(I_k;H)\cap C^0(\bar I_k;V)\cap L^2(I_k;\hth),
\end{eqnarray}
have been constructed for $1\le k\le n<N$, then we look for the pair $(\xi_{n+1},\eta_{n+1})$
that coincides with $(\xi_n,\eta_n)$ in $I_n$, and note that the linear parabolic problem (\ref{eq:4.23}), (\ref{eq:4.24}) has a unique solution $\xi_{n+1}$ on $\Omega\times I_{n+1}$ that satisfies (\ref{eq:4.25}) for $k=n+1$. Inserting
$\xi_{n+1}$ in (\ref{eq:4.21}) 
(where $n$ is replaced by $n+1$), we then find that the linear parabolic problem 
(\ref{eq:4.21}), (\ref{eq:4.22}) admits a
unique solution $\eta_{n+1}$ that fulfills (\ref{eq:4.26}) for $k=n+1$. Hence, we conclude that 
$(\xi^\tau ,\eta^\tau)=(\xi_N,\eta_N)$ satisfies (\ref{eq:4.17})--(\ref{eq:4.20}), and (\ref{eq:4.11}), 
(\ref{eq:4.12}).  

\vspace{3mm}
\underline{Step 2: \,\,A priori estimates.} \,\quad We now prove a series of a priori estimates for the functions $(\xi^\tau,\eta^\tau)$. In the following, we denote by $C_i$ ($i\in\nz$) some generic positive constants, which may depend on $\,\varepsilon, \delta,\rho_*, \rho^*,\mu^*, T, K^*_1, K^*_2,$ $K_3^*$, but not on $\tau$ (i.\,e., not on $N$). For the sake of simplicity,
we omit the superscript $\tau$ and simply write $(\xi,\eta)$.  

\vspace{2mm}
\underline{First a priori estimate.} \,\quad Observe that $\,2\,\bar\rho\,\eta\,\eta_t=
\left(\bar\rho\,\eta^2\right)_t-\bar\rho_t\,\eta^2$. Hence, testing (\ref{eq:4.17}) by $\eta$, and invoking
(\ref{eq:4.19}) and Young's inequality, we have, for 
$0\le t\le T$,  
\begin{eqnarray}
\label{eq:4.27}
&&\int_\Omega \left(\frac {\ve}{2}+\bar\rho(t)\right)\eta(t)^2\,dx \,+\, \int_0^t \|\nabla\eta(s)\|_H^2\,ds
\,+\,\int_0^t\!\!\int_\Gamma \alpha\,\eta^2\,d\sigma\,ds\nonumber\\
&&\le \int_0^t\!\!\int_\Gamma \frac {\alpha}2 \,\eta^2\,d\sigma\,ds\,+\,C_1\int_0^t\!\!\int_\Gamma h^2\,d\sigma\,ds\nonumber\\
&&\quad+\,2\int_0^t\!\!\int_\Omega
|\bar\mu_t|\,|\xi|\,|\eta|\,dx\,ds \,+\int_0^t\int_\Omega
|\bar\mu|\,|\xi_t|\,|\eta|\,dx\,ds \,.
\end{eqnarray}
For any $\gamma>0$, we have, by Young's inequality and (\ref{eq:4.4}), that
\begin{eqnarray}
\label{eq:4.28}
&&\int_0^t\int_\Omega
|\bar\mu|\,|\xi_t|\,|\eta|\,dx\,ds \,\le\,\|\bar\mu\|_{L^\infty(Q)}
\int_0^t\|\eta(s)\|_H\,\|\xi_t(s)\|_H\,ds \nonumber\\
&&\quad\le\,\gamma \int_0^t\|\xi_t(s)\|_H^2\,ds\,+\,\frac {C_2}{\gamma}\int_0^t\|\eta(s)\|_H^2\,ds\,.
\end{eqnarray}
Moreover,
\begin{eqnarray}
\label{eq:4.29}
&&\int_0^t\int_\Omega |\bar\mu_t|\,|\xi|\,|\eta|\,dx\,ds\le \int_0^t\|\bar\mu_t(s)\|_H\,\|\xi(s)\|_
{L^4(\Omega)}\,\|\eta(s)\|_{L^4(\Omega)}\,ds\nonumber\\
&&\quad \le\,\gamma\int_0^t\|\eta(s)\|^2_V\,ds\,+\,\frac{C_3}{\gamma}\,\int_0^t\|\bar\mu_t(s)\|_H^2\,
\|\xi(s)\|^2_V\,ds\,.
\end{eqnarray}
Notice that, by virtue of (\ref{eq:4.5}), the mapping $\,s\mapsto \|\bar\mu_t(s)\|_H^2\,$ is bounded 
by a function in $L^1(0,T)$.

Next, we add $\xi$ on both sides of Eq. (\ref{eq:4.18}) and test the resulting equation by $\xi_t$. On
using Young's inequality again, we obtain:
\begin{eqnarray}
\label{eq:4.30}
&&\frac {\delta} 4 \int_0^t\|\xi_t(s)\|_H^2\,ds \,+\,\frac 1 2 \left(\|\xi(t)\|_H^2\,+\,\|\nabla\xi(t)\|_H^2
\right)\nonumber\\
&&\quad\le \,C_4\,\Bigl(\int_0^t\|\eta(s)\|_H^2\,ds \,+\,\int_0^t\|\xi(s)\|_H^2\,ds\Bigr)\,.
\end{eqnarray} 
Combining the inequalities (\ref{eq:4.27})--(\ref{eq:4.30}), and
choosing $\gamma>0$ sufficiently small, we conclude from Gronwall's lemma that 
\begin{eqnarray}
\label{eq:4.31}
&{\displaystyle \int_0^T\left(\|\xi_t(t)\|_H^2\,+\,\|\eta(t)\|_V^2\right)dt\,
+\,\max_{0\le t\le T}\left(
\|\xi(t)\|_V^2\,+\,\|\eta(t)\|_H^2\right)}&\nonumber\\[1mm]
&{\displaystyle \le
C_5\int_0^T\!\!\int_\Gamma |h|^2\,d\sigma\,dt.}&
\end{eqnarray}
Thanks to (\ref{eq:4.19}), we may also infer (possibly by choosing a larger $C_5$) that
\begin{equation}
\label{eq:4.32}
\|\xi(t)\|_W^2 \,\le\, C_5  \Big(\|\Delta \xi(t)\|_H^2 + 
\int_0^T\!\!\int_\Gamma |h|^2\,d\sigma\,dt  \Big) \quad \mbox{for all }\,t\in\,[0,T]\,.
\end{equation}

\vspace{5mm}
\noindent
\underline{Second a priori estimate.} \,\quad We test (\ref{eq:4.17}) by $\eta_t$ and apply Young's
inequality in order to obtain 
\begin{eqnarray}
\label{eq:4.33}
&\hspace*{-15mm}&\ve\int_0^t\|\eta_t(s)\|_H^2\,ds 
\,+\,\frac  1 2\,\|\nabla\eta(t)\|_H^2\,+\,\int_\Gamma \frac{\alpha} 2\,|\eta(t)|^2\,d\sigma
\,\le\,\int_0^t\!\!\int_\Gamma\alpha\,h\,\eta_t\,d\sigma\,ds\nonumber\\[1mm]
&&\quad +\,
\int_0^t\!\!\int_\Omega\left(2\,|\bar\mu_t|\,|\xi|
+|\bar\mu|\,|\xi_t|+|\bar\rho_t|\,|\eta|\,\right)\,|\eta_t|\,dx\,ds.\quad
\end{eqnarray}
By (\ref{eq:4.4}), we can infer from Young's inequality that
\begin{equation}
\label{eq:4.34}
\int_0^t\int_\Omega |\bar\mu|\,|\xi_t|\,|\eta_t|\,dx\,ds \,\le\,
\gamma\int_0^t\|\eta_t(s)\|_H^2\,ds\,+\,\frac{C_6}{\gamma}\int_0^t\|\xi_t(s)\|_H^2\,ds\,.
\end{equation}
Moreover, by virtue of H\"older's and Young's inequalities,
\begin{eqnarray}
\label{eq:4.35}
&\hspace*{-6mm}&\int_0^t\int_\Omega|\bar\rho_t|\,|\eta|\,|\eta_t|\,dx\,ds
\nonumber\\
&\hspace*{-6mm}&
\le\gamma\int_0^t\|\eta_t(s)\|_H^2\,ds+\frac{C_7}{\gamma}\int_0^t
\|\bar\rho_t(s)\|^2_{L^4(\Omega)}\|\eta(s)\|^2_{L^4(\Omega)}\,ds
\nonumber\\
&\hspace*{-6mm}&\le \gamma\int_0^t\|\eta_t(s)\|_H^2\,ds\,+\,\frac{C_8}{\gamma}\int_0^t\,
\|\bar\rho_t(s)\|^2_V\,\|\eta(s)\|^2_V\,ds \,.
\end{eqnarray}
Observe that by (\ref{eq:4.5}) the mapping $\,s\mapsto \|\bar\rho_t(s)\|^2_V\,$ is
bounded by a function in $L^1(0,T)$.

Also, we have, owing to the continuity of the embedding $W\subset L^\infty(\Omega)$ and
(\ref{eq:4.32}), 
\begin{eqnarray}
\label{eq:4.36}
\hskip-1cm&&\int_0^t\int_\Omega
2\,|\bar\mu_t|\,|\xi|\,|\eta_t|\,dx\,ds\,\nonumber\\[1mm]
\hskip-1cm&&\le\,\gamma\int_0^t\|\eta_t(s)\|_H^2\,ds\,+\,
\,\frac{C_9}{\gamma}\,\int_0^t\,
\|\bar\mu_t(s)\|^2_H\,\|\xi(s)\|_{L^\infty(\Omega)}^2ds\nonumber\\[1mm]
\hskip-1cm&&\le\,\gamma\int_0^t\!\!\|\eta_t(s)\|_H^2\,ds\,+\,\frac{C_{10}}{\gamma}\Bigl(\int_0^T\!\!
\int_\Gamma |h|^2\,d\sigma\,dt\nonumber\\[1mm]
\hskip-1cm&&\hskip5.5cm+\,
\int_0^t\!\!\|\bar\mu_t(s)\|_H^2\,\|\Delta\xi(s)\|_H^2\,
ds\Bigr)\,,
\end{eqnarray}
where, owing to (\ref{eq:4.5}), the mapping $\,s\mapsto\|\bar\mu_t(s)\|_H^2\,$ 
is bounded by a function in $L^1(0,T)$.

Finally, we employ integration by parts, Young's inequality, (\ref{eq:4.31}), and the trace theorem to obtain
\begin{eqnarray}
\label{eq:4.37}
&&\Bigl|\int_0^t\!\!\int_\Gamma \alpha\,h\,\eta_t\,d\sigma\,ds\Bigr|\,\le\,\int_\Gamma\alpha\,|h(t)|\,|\eta(t)|\,d\sigma
\,+\,\int_0^t\!\!\int_\Gamma\alpha\,|h_t|\,|\eta|\,d\sigma\,ds\nonumber\\[1mm]
&&\le\,\int_\Gamma \frac{\alpha}4\,|\eta(t)|^2\,d\sigma\,+\,C_{11}\,\|h\|^2_{\heg}\,.
\end{eqnarray}

Next, we formally test (\ref{eq:4.18}) by $\,-\Delta\xi_t\,$ to obtain, for every $t\in [0,T]$,
\begin{eqnarray}
\label{eq:4.38}
\delta\!\!\int_0^t\!\|\nabla\xi_t(s)\|_H^2\,ds+\frac 1 2 \|\Delta\xi(t)\|_H^2 =\int_0^t\!\!\int_\Omega
\left(-\left(\mathcal{T}_\tau\eta\right)-f''(\bar\rho)\,\xi\right)\,
\Delta\xi_t\,dx\,ds.\quad
\end{eqnarray}

Now, by virtue of (\ref{eq:4.31}) 
and invoking Young's inequality, we have
\begin{eqnarray}
\label{eq:4.39}
&&\Bigl|\int_0^t\!\!\int_\Omega \left(\mathcal{T}_\tau \eta\right) 
\Delta\xi_t\,dx\,ds\Bigr|\nonumber\\[1mm]
&&\le\,
\int_\Omega\!\! \left|\left(\mathcal{T}_\tau \eta\right)(t)\right|\,|\Delta\xi(t)|\,dx\,
+\int_0^t\!\!\int_\Omega \left|\partial_t\left(\mathcal{T}_\tau\eta\right)\right|\,|\Delta\xi|\,dx\,ds
\nonumber\\[1mm]
&&\le\,\frac 1 8\,\|\Delta\xi(t)\|_H^2\,
\,+\,C_{12} \int_0^T\!\!\int_\Gamma |h|^2\,d\sigma\,dt\nonumber\\[1mm]
&&\quad   {}+\,\gamma\int_0^t\!\|\eta_t(s)\|_H^2\,ds\,
+\,\frac 1 {4\gamma} \int_0^t\!\! \|\Delta\xi (s)\|_{H}^2 \,ds\,  .
\end{eqnarray}
Moreover, it turns out that
\begin{eqnarray}
\label{eq:4.40}
\hspace*{-8mm}\Bigl|\int_0^t\!\!\int_\Omega f''(\bar\rho)\,\xi\,\Delta\xi_t\,dx\,ds\Bigr|&\!\!\le\!\!&
\int_\Omega\!\!|f''(\bar\rho(t))|\,|\xi(t)|\,|\Delta\xi(t)|\,dx\nonumber\\[1mm]
&&+\int_0^t\!\!\int_\Omega \left|f'''(\bar\rho)\,\bar\rho_t\,\xi + f''(\bar\rho)\,\xi_t\right|
\,|\Delta\xi|\,dx\,ds.\quad
\end{eqnarray}

We have, owing to (\ref{eq:4.4}) and (\ref{eq:4.31}),
\begin{equation}
\label{eq:4.41}
\int_\Omega\!\!|f''(\bar\rho(t))|\,|\xi(t)|\,|\Delta\xi(t)|\,dx\,\le\,\frac 1 8 \,\|\Delta\xi(t)\|_H^2
\,+\,C_{13}\int_0^T\!\!\int_\Gamma |h|^2\,d\sigma\,dt\,.
\end{equation}

Also the second integral on the right-hand side of (\ref{eq:4.40}) is bounded, 
since (\ref{eq:4.4}), (\ref{eq:4.5}), and (\ref{eq:4.31}) 
imply that 
\begin{eqnarray}
\label{eq:4.42}
\hskip-1cm&&\int_0^t\!\!\int_\Omega \left|f'''(\bar\rho)\,\bar\rho_t\,\xi + f''(\bar\rho)\,\xi_t\right| \,|\Delta\xi|\,dx\,ds \nonumber \\[1mm]
\hskip-12cm&& \leq \, 
C_{14} \,\int_0^t\!\! 
\bigl( \|\bar\rho_t (s)\|_{L^4(\Omega)}^2 \, \|\xi(s)\|_{L^4(\Omega)}^2 
+ \|\xi_t(s)\|_H^2
\bigr) \,ds
+ \int_0^t\!\! \|\Delta\xi (s)\|_{H}^2 \,ds\, 
\nonumber \\[1mm]
\hskip-1cm&& \leq \, 
C_{15} 
\Bigl( \max_{0\le t\le T}\,\|\xi(t)\|_{V}^2 \int_0^t\!\! \|\bar\rho_t (s)\|_{V}^2 \,ds  
+ \int_0^t \|\xi_t(s)\|_H^2 \, ds
\Bigr)+ \int_0^t\!\! \|\Delta\xi (s)\|_{H}^2 \,ds\, 
\nonumber \\[1mm]
\hskip-1cm&&
\leq \, C_{16} \int_0^T\!\!\int_\Gamma |h|^2\,d\sigma\,dt +
\int_0^t\!\! \|\Delta\xi (s)\|_{H}^2 
\,ds\,, 
\end{eqnarray}
thanks to the continuity of the embedding $V\subset L^4(\Omega)$.
Thus, combining the estimates (\ref{eq:4.33})--(\ref{eq:4.42}), choosing $\gamma>0$ sufficiently
small, and invoking Gronwall's inequality, we can infer that
\begin{eqnarray}
\label{eq:4.43}
&{\displaystyle
  \int_0^T\!\!\Bigl(\|\eta_t(t)\|_H^2\,+\,\|\xi_t(t)\|_V^2\Bigr)
\,dt\,+\,
\max_{0\le t\le
  T}\,\left(\|\eta(t)\|_V^2\,+\,\|\xi(t)\|_W^2\right)}&\nonumber\\[1mm]
&{\displaystyle\le\,C_{17}\,\|h\|_{\heg}^2\,.   }
\end{eqnarray}

Next, we compare terms in  (\ref{eq:4.17}) and, arguing as in
the derivation of  (\ref{eq:4.33})--(\ref{eq:4.37}), we readily find that
\begin{equation}
\int_0^T \|\Delta\eta(t)\|_H^2\,dt\,\le\,C_{18}\,\|h\|_{\heg}^2 \,.\nonumber
\end{equation}
Thus, by owing to elliptic regularity (cf.~(\ref{eq:4.19}) and Remark~2.2), we conclude that
\begin{equation}
\label{eq:4.44}
\int_0^T \|\eta(t)\|_{\hth}^2\,dt\,\le\,C_{19}\,\|h\|_{\heg}^2\,.
\end{equation}

Finally, we differentiate Eq. (\ref{eq:4.18}) with respect to~$t$. We obtain:
\begin{equation}
\label{eq:4.45}
\delta\,\xi_{tt}-\Delta\xi_t=\partial_t(\mathcal{T}_\tau\eta)-f'''(\bar\rho)\,
\bar\rho_t\,\xi
-f''(\bar\rho)\,\xi_t\,\quad\mbox{a.\,e. in }\,Q.
\end{equation}

From (\ref{eq:4.4})--(\ref{eq:4.6}), (\ref{eq:4.43}) and (\ref{eq:4.44}), we can infer that we may test (\ref{eq:4.45})
by any of the functions $\,\xi_t$, $-\Delta\xi_t$, and $\,\xi_{tt}$, in order to find that
\begin{equation}
\label{eq:4.46}
\int_0^T\!\!\Bigl(\|\xi_{tt}(t)\|_H^2\,+\,\|\Delta\xi_t(t)\|_H^2\Bigr)\,dt\,+\,
\max_{0\le t\le T}\,\|\xi_t(t)\|_V^2\,
\le\,C_{20}\,\|h\|_{\heg}^2\,.
\end{equation}  

\vspace{2mm}
\underline{Step 3: \,\,Passage to the limit.} \,\quad Let  $(\xi^\tau,\eta^\tau)$ denote the solution
to the system (\ref{eq:4.17})--(\ref{eq:4.20}) associated with $\tau=T/N$, 
for $N\in\nz$. In Step 2, we have shown that
there is some $C>0$, which does not depend on $\tau$, such that
\begin{eqnarray}
\label{eq:4.47}
\hskip-1cm &&\|\xi^{\tau}\|_{H^2(0,T;H)\cap C^1([0,T];V)\cap H^1(0,T;W)\cap C(\overline
  Q)}\nonumber\\[1mm]
\hskip-1cm &&{}+\,
\|\eta^{\tau}\|_{H^1(0,T;H)\cap C^0([0,T];V)\cap L^2(0,T;\hth)}\,\le\,C\,\|h\|_{\heg}\,.
\end{eqnarray}
Hence, there is a subsequence  $\tau_k \searrow 0$ such that 
\begin{eqnarray}
\label{eq:4.48}
\hskip-1cm &&\xi^{\tau_k}\to \xi \ \mbox{ weakly star in }\,\,\,
H^2(0,T;H)\cap W^{1,\infty}(0,T;V)\cap H^1(0,T;W),\nonumber\\[2mm]
\hskip-1cm &&\eta^{\tau_k}\to \eta \ \mbox{ weakly star in }\,\,\, H^1(0,T;H)\cap L^\infty(0,T;V)\cap L^2(0,T;\hth).
\nonumber\\
\hskip-1cm &&
\end{eqnarray}
From the trace theorem we can infer that 
$\xi$ satisfies the boundary 
condition given in (\ref{eq:4.9}),
while the boundary condition for $\eta$ will be satisfied 
(either in the variational sense or in the sense of the appropriate trace theorem,
see Remark~2.2)
once we prove that we can pass to the limit in the products of (\ref{eq:4.7}),
as shown below.
Moreover, it is easily seen that also (\ref{eq:4.10}) is fulfilled.
By compact embedding, we also have, in particular,
\begin{equation}
\label{eq:4.49}
\xi^{\tau_k}\to \xi \,\quad\mbox{strongly in }\,C(\overline{Q}), \qquad \eta^{\tau_k}\to\eta\quad\,\mbox{strongly in }\,
L^2(Q),
\end{equation}
so that $\,\bar\rho\,\eta^{\tau_k}_t\to\bar\rho\,\eta_t\,$ and $\,\bar\mu\,\xi^{\tau_k}_t\to\bar\mu\,\xi_t$,
both weakly in $\,L^2(Q)$, $\,f''(\bar\rho)\,\xi^{\tau_k}\to f''(\bar\rho)\,\xi\,$ strongly in $\,L^2(Q)$,
as well as $\,\bar\mu_t\,\xi^{\tau_k}_t\to \bar\mu_t\,\xi_t\,$ and $\,\bar\rho_t\,\eta^{\tau_k}\to
\bar\rho_t\,\eta$, both strongly in $L^1(Q)$. Finally, it is easily verified that $\,\{\mathcal{T}_{\tau_k}
\eta^{\tau_k}\}\,$ converges strongly in $L^2(Q)$ to $\eta$. In conclusion, we may pass to the limit as
$k \to\infty$ in the system (\ref{eq:4.17})--(\ref{eq:4.20}) (written for $\tau_k $) to find that the pair $(\xi,\eta)$
is in fact a solution to the linearized system (\ref{eq:4.7})--(\ref{eq:4.10}).

We now show the uniqueness. If $(\xi_1,
\eta_1)$, $(\xi_2,\eta_2)$ are two 
solutions having the above properties, then the pair $(\xi,\eta)$, where
$\xi=\xi_1-\xi_2$ and $\eta=\eta_1-\eta_2$, satisfies (\ref{eq:4.7})--(\ref{eq:4.10}) with $h=0$. We thus may repeat
the first a priori estimate in Step 2 to conclude that $\xi=\eta=0$. 

Finally, taking the limit as $\tau \searrow 0$ in (\ref{eq:4.47}) and invoking the lower semicontinuity of norms,
we obtain the inequality (\ref{eq:4.15}). 
This concludes the proof.\linebreak 
\hspace*{2cm} \hfill{\qed}

\subsubsection{Fr\'echet differentiability of the control-to-state mapping}

In this section, we prove the following result.

\vspace{3mm}
{\bf Proposition 4.8} \quad {\em Suppose that the assumptions} (A1)--(A8) {\em are satisfied.
Then the solution operator $S$, viewed
as a mapping from $\mathcal{X}$ to $\mathcal{Y}$, is Fr\'echet differentiable on $\mathcal{U}$.
For any $\bar u\in \mathcal{U}$ the Fr\'echet derivative $DS(\bar u)$ is for $h\in \mathcal{X}$ 
given by $DS(\bar u)h=(\xi,\eta)$, where $(\xi,\eta)$ is the unique solution to the linearized system}
(\ref{eq:4.7})--(\ref{eq:4.10}).

\vspace{3mm}
{\em Proof.} \,\quad Let $\bar u\in \mathcal{U}$ be given and $(\bar \rho,\bar\mu)=S(\bar u)$. Since $\mathcal{U}$ is
an open subset of $\mathcal{X}$, there is some $\lambda>0$ such that $\bar u+h\in \mathcal{U}$ whenever $h\in \mathcal{X}$
satisfies $\|h\|_\mathcal{X}\le\lambda$. In the following, we consider such perturbations $h\in \mathcal{X}$, and we define
$(\rho^h,\mu^h):=S(\bar u+h)$ and put
\begin{equation}
\label{eq:4.50}
z^h:=\mu^h-\bar\mu-\eta^h\,,\quad y^h:=\rho^h-\bar\rho-\xi^h,
\end{equation}
where $(\xi^h,\eta^h)$ denotes the unique solution to the linearized system (\ref{eq:4.7})--(\ref{eq:4.10})
associated with~$h$. 
Since the linear mapping $h\mapsto (\xi^h,\eta^h)$ is by Proposition 4.6 continuous from $\mathcal{X}$ into $\mathcal{Y}$, 
it obviously suffices to show that there is an increasing function 
$g:[0,\lambda]\to [0,+\infty)$ which satisfies $\,\lim_{\,r\searrow 0}\,g(r)/r^2=0\,$
and
\begin{eqnarray}
\label{eq:4.51}
&&\|y^h\|^2_{H^1(0,T;H)\cap C^0([0,T];V)\cap L^2(0,T;W)}\,+\,
\|z^h\|^2_{C^0([0,T];H)\cap L^2(0,T;V)}\nonumber\\[2mm]
&&\le \,g\left(\|h\|_{\heg}\right).
\end{eqnarray}
Using the state system (\ref{eq:1.2})--(\ref{eq:1.4}) and the linearized system 
(\ref{eq:4.7})--(\ref{eq:4.10}), we easily verify
that for $h\in \mathcal{X}$ with $\|h\|_\mathcal{X}\le\lambda$ the pair $(y^h,z^h)$ is a strong solution
to the system
\begin{eqnarray}
\label{eq:4.52}
(\ve+2\bar\rho)\,z^h_t+\bar\rho_t\,z^h+\bar\mu\,y^h_t
+2\bar\mu_t\,y^h-\Delta z^h\hspace{3.5cm}\nonumber\\[1mm]
=
-2\left(\mu^h_t-\bar\mu_t\right)\left(\rho^h-\bar\rho\right)-
\left(\rho^h_t-\bar\rho_t\right)\left(\mu^h-\bar\mu\right)\, \quad
\mbox{a.\,e. in }\,Q,\\[2mm]
\label{eq:4.53}
\delta y^h_t-\Delta y^h+f'(\rho^h)-f'(\bar\rho)- f''(\bar\rho)\,\xi^h
=z^h,\,\quad\mbox{a.\,e. in }\,Q,\\[2mm]
\label{eq:4.54}
\frac{\partial y^h}{\partial {\bf n}} =0,\quad \frac{\partial z^h}{\partial {\bf n}}=-\alpha\,z^h,
\,\quad\mbox{a.\,e. on }\,\Sigma,\\[2mm]
\label{eq:4.55}
y^h(x,0)=z^h(x,0)=0 \,\quad\mbox{for a.\,e. } \, x\in \Omega .
\end{eqnarray}
Notice that 
\begin{eqnarray*}
&&y^h\in H^1(0,T;H)\cap C^0([0,T];V)\cap L^2(0,T;W)\cap C(\bar Q),\\[1mm]
&&z^h\in H^1(0,T;H)\cap C^0([0,T];V)\cap L^2(0,T;\hth).
\end{eqnarray*}
For the sake of a better readability, in the following estimates we omit the superscript $h$ of $y^h$ and~$z^h$. 
Also, we denote by $C_i$ ($i\in\nz$) certain positive
constants that only depend
on  $\,\varepsilon, \delta, \rho_*, \rho^*,\mu^*, T, K^*_1, K^*_2, K^*_3, K^*_4$, but not on~$h$.

We now add $y$ on both sides of Eq. (\ref{eq:4.53}) and test the resulting equation by $y_t$. Using Young's
inequality, we find that for all $t\in [0,T]$ it holds
\begin{eqnarray}
\label{eq:4.56}
\frac {\delta} {2}\int_0^t \|y_t(s)\|_H^2\,ds +\frac 1 2 \left(\|\nabla y(t)\|_H^2+\|y(t)\|_H^2\right)
\,\le\,\frac{2}{\delta}\int_0^t \|z(s)\|_H^2\,ds \qquad\nonumber\\
\,+\,C_1 \int_0^t\|y(s)\|^2_H\,ds\,+\,C_2\int_0^t\|(f'(\rho^h)-f'(\bar\rho)- f''(\bar\rho)
\,\xi^h)(s)\|_H^2\,ds\,.\quad
\end{eqnarray}
In order to handle the third term on the right-hand side of (\ref{eq:4.56}), we note that the stability 
estimate (\ref{eq:3.17}) implies, in particular, that
\begin{equation}
\label{eq:4.57}
\|\rho^h-\bar\rho\|_{L^\infty(Q)}^2\,\le\,K^*_2\,\|h\|^2_{\heg}\,,
\end{equation}
that is, $\rho^h\to\bar\rho$ uniformly on $\overline{Q}$ as $\|h\|_{\heg}\to 0$. Since $f\in C^3(0,1)$,
we can infer from Taylor's theorem and (\ref{eq:4.4}) that
\begin{equation}
\label{eq:4.58}
\left|f'(\rho^h)-f'(\bar\rho)-f''(\bar\rho)\,\xi^h\right|\le
\max_{r_*\le\sigma\le r^*}\,\frac{\left| f'''(\sigma)\right|}2
\left|\rho^h-\bar\rho\right|^2
+|f''(\bar\rho)|\,|y|\,\,\,\,\mbox{on }\,\overline Q.
\end{equation}
 It then follows from the estimates (\ref{eq:3.17}) and  (\ref{eq:4.56})--(\ref{eq:4.57}) that 
\begin{eqnarray}
\label{eq:4.59}
{\displaystyle \frac {\delta} {2}\int_0^t \|y_t(s)\|_H^2\,ds 
+\frac 1 2\, \|y(t)\|_V^2}
&\!\!\le\!\!&{\displaystyle \frac{2}{\delta}\int_0^t \|z(s)\|_H^2\,ds 
\,+\,C_3 \int_0^t\|y(s)\|^2_H\,ds}\nonumber\\[2mm]
&&+\,C_4\,\|h\|_{\heg}^4\,.
\end{eqnarray}

Next, observe that
$\,2\,\bar\rho\,z\,z_t=\left(\bar\rho\,z^2\right)_t-\bar\rho_t\,z^2$.
Therefore, testing (\ref{eq:4.52}) by $z$ yields for every $t\in [0,T]$ that
\begin{eqnarray}
\label{eq:4.60}
&&{\displaystyle \int_\Omega\left(\frac {\ve} 2 +\bar\rho(t)\right)z^2(t)\,dx\,+\,\int_0^t\|\nabla z(s)\|_H^2\,ds
\,+\,\int_0^t\!\!\int_\Gamma \alpha\,|z|^2\,d\sigma\,dt}\nonumber\\[1mm]
&&=\,-\int_0^t\!\!\int_\Omega\left(\bar\mu\,y_t\,+\,2\,\bar\mu_t
\,y\right)z\,dx\,ds
{\displaystyle \,-\,2\int_0^t\!\!\int_\Omega \left(\mu^h_t-
\bar\mu_t\right)\left(\rho^h-\bar\rho\right)\,z
\,dx\,ds}\nonumber\\[1mm]
&&\quad -\int_0^t\!\!\int_\Omega 
\left(\rho^h_t-\bar\rho_t\right)\left(\mu^h-\bar\mu\right)
\,z\,dx\,ds.
\end{eqnarray} 
We estimate the terms on the right-hand side of (\ref{eq:4.60}) individually. At first, 
using (\ref{eq:4.4}) and Young's
inequality, we find that
\begin{equation}
\label{eq:4.61}
\int_0^t\!\!\int_\Omega |\bar\mu|\,|y_t|\,|z|\,dx\,ds\,\le\,\gamma\int_0^t\|y_t(s)\|_H^2\,ds
\,+\,\frac{C_5}{\gamma}\int_0^t\|z(s)\|_H^2\,ds.
\end{equation}
Moreover, using the continuity of the embedding $H^1(\Omega)\subset L^4(\Omega)$, 
as well as H\"older's and Young's inequalities, we have
\begin{eqnarray}
\label{eq:4.62}
 2 \int_0^t\!\!\int_\Omega |\bar\mu_t|\,|y|\,|z|\,dx\,ds\,\le\,2\int_0^t\|\bar\mu_t(s)\|_H
 \,\|z(s)\|_{L^4(\Omega)}\,\|y(s)\|_{L^4(\Omega)}\,ds\nonumber\\
 \le \,\gamma\int_0^t\|z(s)\|_V^2\,ds\,+\,\frac{C_6}{\gamma}\int_0^t\|\bar\mu_t(s)\|^2_H\,\|
 y(s)\|_V^2\,ds\,.
 \end{eqnarray}
 Observe that by (\ref{eq:2.4}) the mapping $\,s\mapsto \|\bar\mu_t(s)\|_H^2\,$ belongs to $L^1(0,T)$.
 
At this point, we can conclude from (\ref{eq:3.17}) and (\ref{eq:4.57}), invoking Young's inequality, that
 \begin{eqnarray}
 \label{eq:4.63}
 &&\int_0^t\!\!\int_\Omega 2\,\left|\mu_t^h-\bar\mu_t\right|\,\left|\rho^h
 -\bar\rho \right|\,|z|\,dx\,ds\nonumber\\
 &&\le\,2\int_0^t \left\|(\mu_t^h-\bar\mu_t)(s)\right\|_H\,
 \left\|(\rho^h  -\bar\rho)(s) \right\|_{L^\infty(\Omega)}\,\|z(s)\|_H\,ds\nonumber\\
&& \le\,C_7\,\left\|\rho^h  -\bar\rho \right\|_{L^\infty(Q)}^2\int_0^t
 \left\|(\mu_t^h-\bar\mu_t)(s)\right\|_H^2\,ds \,+\,\int_0^t\|z(s)\|^2_H\,ds\nonumber\\
 &&\le \int_0^t\|z(s)\|^2_H\,ds\,+\,C_8\,\|h\|_{\heg}^4\,.
 \end{eqnarray}
 
 Finally, we invoke (\ref{eq:3.17}) and H\"older's and Young's inequalities, as well as the continuity of the 
 embedding $H^1(\Omega) \subset L^4(\Omega)$, to obtain that
  \begin{eqnarray}
  \label{eq:4.64}
 &&\int_0^t\!\!\int_\Omega \left|\rho_t^h-\bar\rho_t\right|\,\left|\mu^h
 -\bar\mu \right|\,|z|\,dx\,ds\nonumber\\
 &&\le \,\max_{0\le s\le t}\|z(s)\|_H \int_0^t\left\|(\rho_t^h-\bar\rho_t)(s)\right\|_{L^4(\Omega)}\,\left\|(\mu^h
 -\bar\mu)(s) \right\|_{L^4(\Omega)}\,ds\nonumber\\
 &&\le\,\gamma \,\max_{0\le s\le t}\|z(s)\|_H^2\,+\,\frac{C_9}{\gamma}\int_0^t\left \|
 (\rho_t^h-\bar\rho_t)(s)\right\|_V^2\,ds\,\,
 \int_0^t\left \|
 (\mu^h-\bar\mu)(s)\right\|_V^2\,ds\nonumber\\[2mm]
 &&\le \,\gamma \,\max_{0\le s\le t}\|z(s)\|_H^2\,+\,C_{10}\,\|h\|_{\heg}^4\,.
 \end{eqnarray}
 
 Combining the estimates (\ref{eq:4.59})--(\ref{eq:4.64}), taking the maximum with respect to $t\in [0,T]$,
 adjusting $\gamma>0$ appropriately small, and invoking Gronwall's lemma, we arrive at the conclusion that $(y^h,z^h)
 =(y,z)$ satisfies the inequality
 \begin{eqnarray}
 \label{eq:4.65} \displaystyle 
 \|y^h\|_{H^1(0,T;H)\cap C^0([0,T];V)}^2\,+\,\|z^h\|^2_{C^0([0,T];H)\cap L^2([0,T];V)} \nonumber\\[2mm]
 \displaystyle \,\le\,C_{11}\,\|h\|_{\heg}^4\,.
 \end{eqnarray}
 Finally, testing (\ref{eq:4.53}) by $\,-\Delta y^h$, and using (\ref{eq:4.58}), we find that also
 \begin{equation}
 \label{eq:4.66}
 \|y^h\|^2_{L^2(0,T;W)}\,\le\, C_{12}\,\|h\|_{\heg}^4\,.
 \end{equation}
 
 Therefore, the function $g(r):=\left(C_{11}+C_{12}\right)\,r^4$ has the requested properties. This concludes the proof of the assertion. \qed
 
 \vspace{7mm}
 {\bf Corollary 4.9} \,\quad {\em Let the assumptions} (A1)--(A8) {\em be fulfilled, and
 let $\bar u\in \uad$  be an optimal control for the problem} {\bf (CP)} {\em with associated
 state $(\bar\rho,\bar\mu)=S(\bar u)$. Then, for every $v\in \uad$, 
 \begin{equation}
 \label{eq:4.67}
 \int_0^T\!\!\!\int_ \Gamma\!\! \beta_1\,\bar u(v-\bar u)\,d\sigma\,dt +\int_\Omega\!\! (\bar\rho(T)-\rho_T)\,
 \xi(T)\,dx\,+\,\int_0^T\!\!\int_\Omega\!\!\beta_2\,(\bar\mu-\mu_T)\,\eta\,dx\,dt\,\ge\,0,\quad
  \end{equation} 
 where $(\xi,\eta)$ is the unique solution to the linearized system} 
 (\ref{eq:4.7})--(\ref{eq:4.10}) {\em associated with 
 $h=v-\bar u$.} 
 
 \vspace{2mm}
 {\em Proof.} \,\quad Let $\,v\in \uad$ be arbitrary and $h=v-\bar u$. Then $\bar u +\lambda h  \in \uad$ for $0<\lambda\le
 1$. For any such $\lambda$, we have
 \begin{eqnarray*}
0&\le&\frac{J(\bar u+\lambda h,S(\bar u+\lambda h))-J(\bar u,S(\bar u))}\lambda\nonumber\\[3mm]
&\le& \frac{J(\bar u+\lambda h,S(\bar u+\lambda h))
-J(\bar u,S(\bar u+\lambda h))}{\lambda}\\[3mm]
 &&   {}+
 \frac{ J(\bar u,S(\bar u+\lambda h))-J(\bar u,S(\bar u))}{\lambda}\,.\qquad
 \end{eqnarray*}
 It follows immediately from the definition of the cost functional $J$ that the first summand on the 
 right-hand side of this inequality converges to $\,\int_0^T\!\!\int_\Gamma \beta_1\,\bar u\,(v-\bar u)\,d\sigma\,dt\,$
 as $\lambda\searrow 0$. For the second summand, we obtain from Proposition 4.8 that
 \begin{eqnarray*}
&&{\displaystyle \lim_{\lambda\searrow 0}\,\frac{ J(\bar u,S(\bar u+\lambda h))-J(\bar 
u,S(\bar u))}{\lambda}}\\
&&{\displaystyle  = \int_\Omega\!\! (\bar\rho(x,T)-\rho_T(x))\,
 \xi(x,T)\,dx\,+\,\int_0^T\!\!\int_\Omega\!\!\beta_2\,(\bar\mu-\mu_T)\,
\eta\,dx\,dt\,,}
 \end{eqnarray*}
 whence the assertion follows. \qed
 
 \subsubsection{The optimality system}
 
 Let $\bar u\in \uad$ be an optimal control for {\bf (CP)} with associated state 
 $(\bar\rho,\bar\mu)=S(\bar u)$. Then, for every $\,v\in\uad$, (\ref{eq:4.67}) holds. We 
 now aim to eliminate $(\xi,\eta)$ by introducing the adjoint
 state variables. To this end, we consider the {\em adjoint system}\,:
 
 \begin{eqnarray}
 \label{eq:4.68}
 -(\ve+2\bar\rho)\,q_t-\bar\rho_t\,q-\Delta q=p+\beta_2\left(\bar\mu-\mu_T\right)\,\quad
 \mbox{a.\,e. in }\,Q,\\[1mm]
 \label{eq:4.69}
 \frac{\partial q}{\partial \bf n}=-\alpha\,q\quad \mbox{a.\,e. in }\,\Sigma,\qquad q(x,T)=0\, \quad
 \mbox{for a.\,e. } \, x\in \Omega,\\[2mm]
 \label{eq:4.70}
 -\delta p_t-\Delta p+f''(\bar\rho)\,p=\bar\mu\,q_t-\bar\mu_t\,q\,\quad
 \mbox{in }\,Q,\\[1mm]
 \label{eq:4.71}
 \frac{\partial p}{\partial \bf n}=0\quad \mbox{on }\,\Sigma,\qquad  
 \delta\,p(T)=\bar\rho(T)-\rho_T \,\quad \mbox{in }\,\Omega\,, 
 \end{eqnarray}
 
 which is a linear backward-in-time parabolic system for the adjoint state variables $p$ and $q$.
 
 It must be expected that the adjoint state variables $(p,q)$ be less regular than the state
 variables $(\bar\rho,\bar\mu)$. Indeed, we only have $p(T)\in L^2(\Omega)$, and thus 
 (\ref{eq:4.70}) and (\ref{eq:4.71})
 should be interpreted in the ususal weak sense. That is, we look for a vector-valued function
 $\,p\in H^1(0,T;V^*)\cap C^0([0,T];H)\cap L^2(0,T;V)$ that, in addition to the final time condition
 (\ref{eq:4.71}), satisfies 
 \begin{eqnarray}
 \label{eq:4.72}
 &&\langle -\delta\,p_t(t),v\rangle_{V^*,V}\,+\,\int_\Omega\nabla p(t)\cdot\nabla v\,dx \,+\,
 \int_\Omega f''(\bar\rho(t))\,p(t)\,v\,dx\nonumber\\[1mm]
 &&\quad= \int_\Omega \left(\bar\mu(t)\,q_t(t)-\bar\mu_t(t)\,q(t)\right)\,v\,dx\,,
 \end{eqnarray}
 for every $v\in V$ and almost every $t\in (0,T)$. 
  Notice that if $q\in H^1(0,T;H)\cap C^0([0,T];V)$, then it is easily seen that 
 $\,\bar\mu\,q_t-\bar\mu_t\,q\in L^{3/2}(Q)$, so that the integral on the right-hand side of
 (\ref{eq:4.72}) makes sense. On the other hand, if $p$ has the expected regularity then
 the solution to (\ref{eq:4.68}), (\ref{eq:4.69}) should belong to 
 $H^1(0,T;H)\cap C^0([0,T];V)\cap L^2(0,T;\hzwei)$. 
 
 The following result is an analogue of Theorem~3.7 in \cite{CGPS5}.
  
\vspace{7mm}
 {\bf Theorem 4.10} \,\quad {\em Suppose that $\bar u\in\uad$ is an optimal control for} {\bf (CP)}
 {\em with associated state $(\bar\rho,\bar\mu)=S(\bar u)$. Then the adjoint system} 
 (\ref{eq:4.68})--(\ref{eq:4.71}) {\em
 has a unique weak solution $(p,q)$ with}
 $p\in H^1(0,T; V^*)\cap C^0([0,T];H)\cap L^2(0,T;V)$, $q\in
 H^1(0,T;H)\cap C^0([0,T];V)\cap L^2(0,T;\hzwei)${\em ;
 moreover, for any $v\in \uad$, we have the inequality}
 \begin{equation}
 \label{eq:4.73}
 \int_0^T\!\!\int_\Gamma \beta_1\,\bar u\,(v-\bar u)\,d\sigma\,dt\,
+\,\int_0^T \!\!\int_\Gamma \alpha \, q\,(v-\bar u) \,d\sigma\,dt\,\ge\,0\,.
 \end{equation} 
 
 \vspace{2mm}
 {\em Proof.} \,\quad The existence and uniqueness result for the adjoint state variables $p$ and $q$ 
follows using the same line of arguments as in the proof of Proposition~3.6 in \cite{CGPS5}, 
with only minor and straightforward changes that are due to the different boundary condition for $q$. 
Now let $v\in\uad$ be given. A standard calculation 
(which can be left as an easy  exercise to the reader), 
using the linearized system (\ref{eq:4.7})--(\ref{eq:4.10}) with $h=v-\bar u$,  
repeated integration by parts, and the well-known integration by parts formula
 $$\int_0^T\bigl(\langle v_t(t),w(t)\rangle_{V^*,V}\,+\, \langle w_t(t),v(t)\rangle_{V^*,V}\bigr)\,dt
 =\int_\Omega\bigl(v(T)w(T)-v(0)w(0)\bigr)\,dx$$
 (which holds for all functions $v,w\in H^1(0,T;V^*)\cap L^2(0,T;V)$), yields the identity 
 \begin{eqnarray}
 \label{eq:4.74}
 &&\int_\Omega (\bar\rho(x,T)-\rho_T(x))\,\xi(x,T)\,dx\,+\,\int_0^T\!\!\int_\Omega \beta_2(\bar\mu-\mu_T)\eta\,dx\,dt
 \nonumber\\[1mm]
 &&=\,\int_0^T\!\!\int_\Gamma \alpha \, q\,(v-\bar u)\,d\sigma\,dt\,.
 \end{eqnarray}
 
 The variational inequality (\ref{eq:4.73}) is thus a direct consequence of Corollary 4.9. 
 \linebreak \hspace*{1cm} \hfill{\qed}

\vspace{5mm}
\noindent
{\bf Acknowledgement. }
P. Colli and G. Gilardi 
gratefully acknowledge the financial support of
 the MIUR-PRIN Grant 2008ZKHAHN {\em ``Phase transitions, 
hysteresis and multiscaling''} and of the IMATI of CNR in Pavia. 
The work of J. Sprekels was supported by the DFG Research Center 
{\sc Matheon} in
Berlin.

 
\end{document}